\documentclass[11pt,a4paper]{article}
\usepackage[T1]{fontenc}
\usepackage{lmodern}
\usepackage[a4paper,margin=27mm,headheight=14pt]{geometry}
\usepackage{amsmath,amssymb,amsthm,mathtools,mathrsfs,booktabs,array}
\usepackage{enumitem,fancyhdr,xcolor,hyperref}
\usepackage{tikz}
\usetikzlibrary{arrows.meta}
\hypersetup{colorlinks=true,linkcolor=blue!45!black,citecolor=blue!45!black,urlcolor=blue!45!black,
 pdftitle={Restrictions of planar doubling measures to piecewise C1 Jordan domains},
 pdfauthor={Xiaohua Wang, Fan Wen}}
\newtheorem{theorem}{Theorem}[section]
\newtheorem{lemma}[theorem]{Lemma}
\newtheorem{proposition}[theorem]{Proposition}
\newtheorem{corollary}[theorem]{Corollary}
\theoremstyle{definition}
\newtheorem{definition}[theorem]{Definition}
\newtheorem{remark}[theorem]{Remark}
\newtheorem{example}[theorem]{Example}
\newcommand{\R}{\mathbb R}
\newcommand{\Leb}{\mathcal L}
\newcommand{\D}{\mathcal D}
\newcommand{\dist}{\operatorname{dist}}

\newcommand{\restr}{\mathbin{\vrule height 1.4ex depth -0.3ex width 0.07ex\vrule height 0.07ex depth 0ex width 0.65ex}}
\DeclareMathOperator{\Lip}{Lip}
\numberwithin{equation}{section}
\setlist[enumerate]{nosep,leftmargin=2em}
\title{Restrictions of planar doubling measures\\to piecewise $C^1$ Jordan domains\thanks{Supported by the National Natural Science Foundation of China (NSFC), Grant Nos.~12571094 and 1250010257.}}
\author{\normalsize Xiaohua Wang\textsuperscript{1}\quad Fan Wen\textsuperscript{2}\thanks{Corresponding author. Email: \href{mailto:wenfan@jnu.edu.cn}{\nolinkurl{wenfan@jnu.edu.cn}}.}\\
\small \textsuperscript{1}\,China Economics and Management Academy\\
\small Central University of Finance and Economics, China\\
\small \href{mailto:fridman_wxh@163.com}{\nolinkurl{fridman_wxh@163.com}}\\[6pt]
\small \textsuperscript{2}\,Department of Mathematics, Jinan University\\
\small Guangzhou 510632, China\\
\small \href{mailto:wenfan@jnu.edu.cn}{\nolinkurl{wenfan@jnu.edu.cn}}}
\date{}
\begin{document}
\maketitle
\begin{abstract}
We compare the class $\mathscr R(\Omega)$ of restrictions of full-support planar doubling measures with the class $\D(\Omega)$ of internal doubling measures on a piecewise $C^1$ Jordan domain. These classes coincide in the absence of zero-interior-angle cusps. At any such cusp, $\mathscr R(\Omega)\setminus\D(\Omega)$ is nonempty. For power-comparable widths, $\D(\Omega)\subsetneq\mathscr R(\Omega)$; for infinitely flat widths the two classes are disjoint. Under smooth endpoint regularity, finite versus infinite contact order gives these two alternatives. Under only $C^1$ regularity, finite logarithmic order does not determine whether the classes intersect, and examples with nonempty intersection and neither class containing the other also occur. For an individual internal measure, membership in $\mathscr R(\Omega)$ is characterized by a potential inequality for its axial projection. We prove this criterion by a positive-kernel construction and Whitney mass replacement, and also obtain sufficient conditions for variable-power widths and an independent finite-system criterion.
\end{abstract}
\noindent\textbf{Keywords:} doubling measure; restriction; Jordan domain; zero interior angle; cusp width; finite order.\\
\noindent\textbf{MSC[2020]:} 28A75

\section{Introduction and main results}

Doubling measures are a basic tool in analysis on metric spaces; see Heinonen~\cite{Hei}. The existence results of Vol'berg and Konyagin~\cite{VK} and Luukkainen and Saksman~\cite{LS}, and the nested-cube construction of K\"aenm\"aki, Rajala, and Suomala~\cite{KRS}, ensure that compact subsets of Euclidean space carry doubling measures. They do not imply that restricting a given ambient doubling measure preserves doubling. We study this question together with the relation between internal doubling measures and planar restrictions, using two classes of measures on the same domain.

All distances are Euclidean. For $X\subset\R^2$, write $B_X(x,r)=X\cap B(x,r)$ and denote by $\D(X)$ the full-support doubling Radon measures on $X$. Thus $\mu\in\D(X)$ means that some $D_\mu\ge1$ satisfies
\begin{equation}\label{eq:doubling}
0<\mu(B_X(x,2r))\le D_\mu\mu(B_X(x,r))<\infty
\qquad(x\in X,\ r>0).
\end{equation}
For a bounded Jordan domain $\Omega$, put $E=\overline\Omega$ and define
\begin{equation}\label{eq:goodclass}
\mathscr R(\Omega)=\{\nu\restr\Omega:\nu\in\D(\R^2)\}.
\end{equation}
The same notation is used for closed domains and truncated models. Both $\mathscr R(\Omega)$ and $\D(\Omega)$ consist of measures on $\Omega$; membership in $\mathscr R(\Omega)$ does not require internal doubling. We study their intersection and the two differences. Thus $\mathscr R(\Omega)\subseteq\D(\Omega)$ means that every planar doubling restriction is doubling, whereas $\D(\Omega)\subseteq\mathscr R(\Omega)$ means that every internal doubling measure is a planar restriction. Equality of measures always means equality on every Borel set.

\subsection{The boundary convention and orders of contact}

\begin{definition}\label{def:domain}
A \emph{piecewise $C^1$ Jordan domain} is a bounded Jordan domain $\Omega$ whose boundary is a cyclic union of finitely many regular $C^1$ arcs up to their endpoints. Each arc admits an injective $C^1$ parametrization on a closed interval with nonvanishing derivative throughout that interval, and distinct arcs meet only at the endpoints of adjacent arcs. A junction $p\in\partial\Omega$ is a \emph{cusp of zero interior angle} if the two tangent rays directed away from $p$ coincide and $\Omega$ lies locally between the incident arcs.
\end{definition}

We first record why the choice between the open domain $\Omega$ and its closure $E=\overline\Omega$ is immaterial. Every $\mu\in\D(E)$ and every $\nu\in\D(\R^2)$ gives zero mass to $\partial\Omega$. To see this, let $m$ be a doubling measure on either connected space $X=E$ or $X=\R^2$. For each $x\in X$ and sufficiently small $r$, connectedness provides $q\in X$ with $|q-x|=3r/4$. Since $B_X(q,r/8)$ lies in $B_X(x,r)\setminus B_X(x,r/2)$ and $B_X(x,r)\subset B_X(q,2r)$, four applications of doubling yield
\begin{equation}\label{eq:nonatomiccontraction}
 m(B_X(x,r/2))\le(1-D_m^{-4})m(B_X(x,r)).
\end{equation}
Iteration gives $m(\{x\})=0$. Away from the finitely many arc junctions, a boundary point has a $C^1$ graph neighborhood. Each sufficiently small ball $B(x,r)$ centered there contains an interior ball $B(a,\kappa r)\subset\Omega\cap B(x,r)$, with a fixed local $\kappa>0$. Doubling assigns this interior ball a fixed positive fraction of $m(B_X(x,r))$. The boundary therefore cannot have density one at such a point, and the differentiation theorem for doubling measures~\cite{Hei} makes the regular boundary arcs null. The finitely many junctions are null by nonatomicity.

Consequently, restriction gives a bijection from $\D(E)$ onto $\D(\Omega)$. Indeed, its inverse assigns zero mass to the boundary: for $\lambda\in\D(\Omega)$, put $\bar\lambda(A)=\lambda(A\cap\Omega)$. This is a finite Radon measure because $\Omega$ is bounded and $\lambda$ is finite on every relative ball. Given $x\in E$ and $r>0$, choose $a\in\Omega$ with $|a-x|<r/4$. Then $B(x,2r)\cap\Omega\subset B_\Omega(a,3r)$ and $B_\Omega(a,r/2)\subset B(x,r)\cap\Omega$, so $\bar\lambda$ is doubling with constant at most $D_\lambda^3$. Ambient boundary nullity also gives
\begin{equation}\label{eq:openclosed}
 \nu\restr\Omega\in\D(\Omega)
 \quad\Longleftrightarrow\quad
 \nu\restr E\in\D(E),\qquad \nu\in\D(\R^2).
\end{equation}
Under this identification, both $\mathscr R(\Omega)$ and $\D(\Omega)$ have the same meaning on the open and closed domains. We use $\Omega$ in the statements and $E$ when compactness is convenient. Since $E$ is compact and geometrically doubling, the existence theorem of Luukkainen and Saksman~\cite{LS} also gives $\D(\Omega)\ne\varnothing$.

A cusp in the complement has interior angle $2\pi$ and is not a zero-interior-angle cusp of $\Omega$. A globally regular $C^1$ Jordan curve has no zero-angle cusp.

If a piecewise $C^1$ Jordan domain has no zero-interior-angle cusp, then $\mathscr R(\Omega)=\D(\Omega)$. Indeed, regular boundary points and nonzero-angle junctions, including those of angle $2\pi$, satisfy a uniform small-scale interior ball condition. Lemma~\ref{lem:cork} gives $\mathscr R(\Omega)\subseteq\D(\Omega)$, while the reflected mass construction of Lemma~\ref{lem:corkrealization} gives the reverse inclusion. The boundary identification above applies to both.

Throughout the main results, unless stated otherwise, $\Omega$ is a piecewise $C^1$ Jordan domain whose unique zero-interior-angle cusp is at $0$, called a \emph{single-cusp domain}.

\begin{figure}[htbp]
\centering
\begin{tikzpicture}[x=1cm,y=1cm,font=\small,
  boundary/.style={line width=.85pt},
  axis/.style={-{Stealth[length=1.8mm]},line width=.45pt},
  measure/.style={{Stealth[length=1.5mm]}-{Stealth[length=1.5mm]},line width=.55pt}]
  \begin{scope}
    \fill[black!7]
      plot[domain=0:3.25,samples=70] (\x,{.20*\x*\x}) --
      plot[domain=3.25:0,samples=70] (\x,{-.12*\x*\x}) -- cycle;
    \draw[axis] (-.25,0) -- (4.15,0) node[right] {$t$};
    \draw[axis] (0,-1.65) -- (0,2.45) node[above] {$y$};
    \draw[boundary] plot[domain=0:3.25,samples=70] (\x,{.20*\x*\x});
    \draw[boundary] plot[domain=0:3.25,samples=70] (\x,{-.12*\x*\x});
    \draw[dashed,line width=.6pt]
      plot[domain=0:3.25,samples=70] (\x,{.04*\x*\x});
    \node[right] at (3.25,2.1125) {$g_+(t)$};
    \node[right] at (3.25,-1.2675) {$g_-(t)$};
    \node[right] at (3.25,.4225) {$h(t)$};
    \node at (2.55,.82) {$\Omega$};
    \fill (0,0) circle[radius=1pt];
    \node[below left] at (0,0) {$0$};
    \node at (1.85,-2.15) {(a) Original cusp};
  \end{scope}
  \draw[axis,line width=.65pt] (5.05,.65) -- (6.35,.65)
    node[midway,above=3pt] {$T$};
  \begin{scope}[xshift=7.15cm]
    \fill[black!7]
      plot[domain=0:3.25,samples=70] (\x,{.16*\x*\x}) --
      plot[domain=3.25:0,samples=70] (\x,{-.16*\x*\x}) -- cycle;
    \draw[axis] (-.25,0) -- (4.15,0) node[right] {$t$};
    \draw[axis] (0,-1.95) -- (0,2.45) node[above] {$y$};
    \draw[boundary] plot[domain=0:3.25,samples=70] (\x,{.16*\x*\x});
    \draw[boundary] plot[domain=0:3.25,samples=70] (\x,{-.16*\x*\x});
    \node[right] at (3.25,1.69) {$f(t)$};
    \node[right] at (3.25,-1.69) {$-f(t)$};
    \draw[densely dotted] (2.3,-.8464) -- (2.3,0);
    \draw[measure] (2.3,0) -- (2.3,.8464)
      node[midway,right=3pt] {$f(t_0)$};
    \node[below right] at (2.3,0) {$t_0$};
    \node at (2.85,-.90) {$\Omega_f$};
    \fill (0,0) circle[radius=1pt];
    \node[below left] at (0,0) {$0$};
    \node at (1.85,-2.15) {(b) Straightened cusp};
  \end{scope}
\end{tikzpicture}
\caption{Local geometry of a zero-interior-angle cusp. Left: the boundary branches
$g_-$ and $g_+$ and their midpoint graph $h$. Right: the shear
$T(t,y)=(t,y-h(t))$ gives $\Omega_f$. The common tangent ray at $0$ is the positive
$t$-axis; the total vertical width at $t_0$ is $2f(t_0)$. Only a cusp neighborhood
is shown. The profiles are schematic; monotonicity is not assumed.}
\label{fig:cusp}
\end{figure}
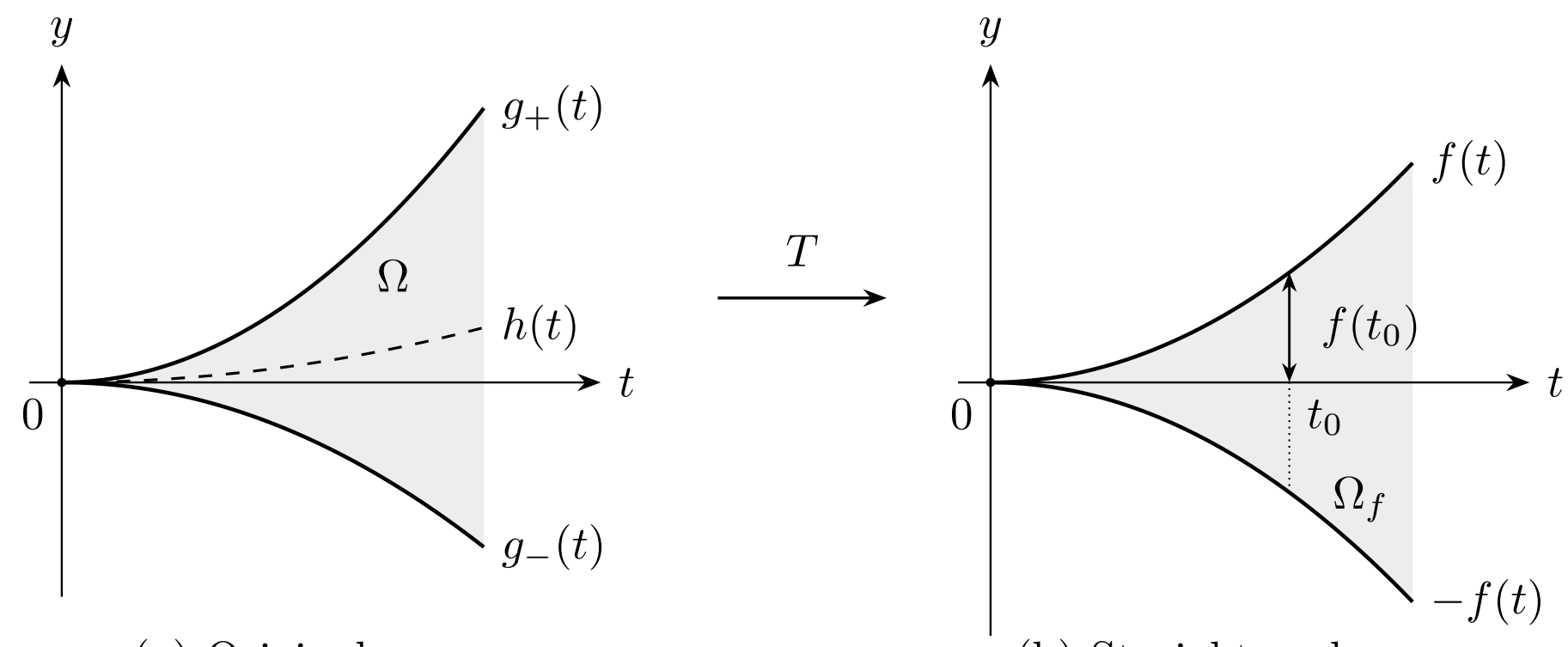

Choose the common tangent as the positive $t$-axis. Near $0$ the boundary arcs are graphs $g_-(t)<g_+(t)$, where $g_\pm\in C^1([0,a])$ and $g_\pm(0)=g'_\pm(0)=0$. Set
\begin{equation}\label{eq:width}
h(t)=\frac{g_+(t)+g_-(t)}2,\qquad
f(t)=\frac{g_+(t)-g_-(t)}2.
\end{equation}
A Lipschitz shear straightens the midpoint graph. The resulting local domain is
\begin{equation}\label{eq:model}
\Omega_f=\{(t,y):0<t<a,\ |y|<f(t)\},\qquad
f(0)=f'(0)=0,\quad f(t)>0\ (t>0).
\end{equation}
The notation $\Omega_f$ describes only the cusp neighborhood; a regular cap completes the bounded domain away from $0$.
Figure~\ref{fig:cusp} illustrates the original cusp and its straightened local model.

We distinguish three conditions on $f$. It has \emph{power-comparable order} $\alpha>1$ if, near $0$,
\begin{equation}\label{eq:powercomparison}
c t^\alpha\le f(t)\le C t^\alpha
\end{equation}
for some $c,C>0$. It has \emph{logarithmic order} $\alpha$ if
\begin{equation}\label{eq:logorder}
\lim_{t\downarrow0}\frac{\log f(t)}{\log t}=\alpha.
\end{equation}
It is \emph{infinitely flat} if, for every $M>0$, there are $C_M,a_M>0$ such that
\begin{equation}\label{eq:flat}
f(t)\le C_Mt^M\qquad(0<t<a_M).
\end{equation}
Power comparison implies the corresponding logarithmic order, but the converse fails. These definitions use the fixed tangent coordinate in \eqref{eq:width}; no derivatives of order greater than one are assumed.

If the incident arcs are $C^\infty$ up to $0$, their \emph{contact order} is
\begin{equation}\label{eq:contact}
m=\min\{j\ge2:f^{(j)}(0)\ne0\},
\end{equation}
with $m=\infty$ if all derivatives vanish. For finite $m$, Taylor's theorem gives $f(t)=b t^m+o(t^m)$ with $b>0$; for $m=\infty$, \eqref{eq:flat} holds.

\subsection{Relations between the two classes}

\begin{theorem}\label{thm:bad}
For every single-cusp piecewise $C^1$ Jordan domain,
$\mathscr R(\Omega)\setminus\D(\Omega)\ne\varnothing$.
A measure in this difference can be chosen as the restriction of an absolutely continuous full-support planar doubling measure.
\end{theorem}

\begin{theorem}\label{thm:flat}
Put $F(r)=\sup_{0\le t\le r}f(t)$. If $\mathscr R(\Omega)\cap\D(\Omega)\ne\varnothing$, then some $c>0$, $N>1$, and $r_0>0$ satisfy
\begin{equation}\label{eq:envelopebound}
F(r)\ge cr^N\qquad(0<r<r_0).
\end{equation}
Consequently, an infinitely flat width gives $\mathscr R(\Omega)\cap\D(\Omega)=\varnothing$.
If $f$ is nondecreasing, $F(r)$ in \eqref{eq:envelopebound} may be replaced by $f(r)$.
\end{theorem}

\begin{theorem}\label{thm:power}
If the width satisfies \eqref{eq:powercomparison} for some $\alpha>1$, then
\[
\D(\Omega)\subsetneq\mathscr R(\Omega),
\qquad \Leb^2\restr\Omega\in\D(\Omega).
\]
Only $C^1$ regularity of the incident arcs is required. In particular, the conclusion holds if $f(t)=bt^\alpha+o(t^\alpha)$ with $b>0$, without any derivative asymptotic assumption.
\end{theorem}

\begin{theorem}\label{thm:smooth}
Suppose the two incident arcs are regular $C^\infty$ up to $0$, and let $m$ be their contact order. If $m<\infty$, then $\D(\Omega)\subsetneq\mathscr R(\Omega)$ and restricted area is doubling. If $m=\infty$, then $\mathscr R(\Omega)\cap\D(\Omega)=\varnothing$.
\end{theorem}

\begin{theorem}\label{thm:counter}
For every $\alpha>1$ there exists a single-cusp piecewise $C^1$ Jordan domain with strictly increasing width
$f\in C^1([0,a])\cap C^\infty((0,a])$, $f(0)=f'(0)=0$, such that
\begin{equation}\label{eq:counterbounds}
t^\alpha\le f(t)\le3t^\alpha\exp\!\sqrt{\log(1/t)}
\quad(0<t<a),
\end{equation}
but $\mathscr R(\Omega)\cap\D(\Omega)=\varnothing$.
Its logarithmic order is $\alpha$.
\end{theorem}

Theorems~\ref{thm:power} and~\ref{thm:counter} show that finite logarithmic order does not determine whether the two classes intersect, even for strictly increasing widths with a polynomial lower bound.

\subsection{Membership of an internal measure}

For the local criterion, truncate and linearly rescale the straightened cusp to
\begin{equation}\label{eq:horn}
\begin{gathered}
E_f=\{(t,y):0\le t\le1,\ |y|\le f(t)\},\\
f\in C^1([0,1]),\quad f(0)=f'(0)=0,\quad f(t)>0\ (t>0).
\end{gathered}
\end{equation}
Other positive axial endpoints are allowed by rescaling. We use
$\pi(t,y)=t$, $I(t,r)=(t-r,t+r)\cap[0,1]$, and the axial projection
\begin{equation}\label{eq:localprojection}
\lambda=\pi_\#\mu\qquad(\mu\in\D(E_f)).
\end{equation}
Restriction of an internal measure to a sufficiently short truncated cusp remains doubling: at the new cross-section, an interior ball of comparable radius gives the required mass comparison; elsewhere small balls inherit the original estimate, and compactness handles the remaining scales. The boundary correspondence following Definition~\ref{def:domain} permits this closed local model while keeping $\Omega$ as the principal object.

\begin{theorem}\label{thm:potentialcriterion}
Let $E_f$ be given by \eqref{eq:horn}, let $\mu\in\D(E_f)$, and set $\lambda=\pi_\#\mu$. Then $\mu\in\mathscr R(E_f)$ if and only if there exist $0<\varepsilon<1$, $N\ge2$, and $C<\infty$ such that, for every sufficiently small $t>0$,
\begin{equation}\label{eq:potentialtesting}
\int_{[0,1]}
\left(\frac{f(t)}{f(t)+|t-s|}\right)^N
f(s)^{-\varepsilon}\,d\lambda(s)
\le C f(t)^{-\varepsilon}\lambda(I(t,f(t))).
\end{equation}
The projection $\lambda$ is doubling and has no atom at $0$. No monotonicity or doubling condition on $f$ is assumed. When the test holds, there is a measure $\nu\in\D(\R^2)$ such that $\nu\restr E_f=\mu$ and $\nu\restr(\R^2\setminus E_f)\ll\Leb^2$.
\end{theorem}

The parameters may depend on the given measure. In particular, this criterion does not assert that the same parameters work for all measures on one domain.

\begin{corollary}\label{thm:axialreduction}
For $\lambda\in\D([0,1])$ define
\begin{equation}\label{eq:canonicalmeasure}
d\eta_\lambda(t,y)=
\frac{\mathbf1_{\{|y|\le f(t)\}}}{2f(t)}\,d\lambda(t)\,dy.
\end{equation}
Then $\eta_\lambda\in\D(E_f)$ and $\pi_\#\eta_\lambda=\lambda$. For any $\mu\in\D(E_f)$,
\[
\mu\in\mathscr R(E_f)\quad\Longleftrightarrow\quad
\eta_{\pi_\#\mu}\in\mathscr R(E_f).
\]
More generally, two internal doubling measures with comparable axial projections either both belong to $\mathscr R(E_f)$ or both lie outside it.
\end{corollary}

\begin{corollary}\label{thm:internalnonrestriction}
If $\mu\in\D(\Omega)\cap\mathscr R(\Omega)$ and $\lambda$ is its local axial projection on $(0,a)$, then, for some $\varepsilon>0$,
\begin{equation}\label{eq:negativewidthmoment}
\int_{(0,a)}f(t)^{-\varepsilon}\,d\lambda(t)<\infty.
\end{equation}
\end{corollary}

This follows already from the potential criterion, since its kernel at any fixed positive $t$ is bounded below on the compact axial interval. A direct proof by strip decay is given in Section~3. The condition is not sufficient; see Remark~\ref{rem:momentnonsufficient}.

\begin{theorem}\label{thm:partialc1}
There exists a single-cusp piecewise $C^1$ Jordan domain whose width belongs to $C^1([0,a])\cap C^\infty((0,a])$ and for which all three sets
\[
\mathscr R(\Omega)\cap\D(\Omega),\qquad
\mathscr R(\Omega)\setminus\D(\Omega),\qquad
\D(\Omega)\setminus\mathscr R(\Omega)
\]
are nonempty.
\end{theorem}

The construction uses narrow width valleys around rescaled compact sets of small upper box dimension. One axial doubling measure passes the potential test; another has no finite negative width moment. Thus nonempty intersection need not imply either inclusion. Proposition~\ref{prop:widthtest} and Theorem~\ref{thm:variablepower} give further sufficient geometric conditions for $\D(\Omega)\subseteq\mathscr R(\Omega)$. Theorem~\ref{thm:dyadiccriterion} gives an independent membership criterion through finite systems of mass constraints.

\subsection{Related work and proof strategy}

Restriction questions belong to the study of ambient doubling measures; see Kaufman and Wu~\cite{KW}. Without boundary control, existence on the closure need not imply existence on the domain: Saksman~\cite{Sak} constructed bounded Jordan domains carrying no doubling measure.

For uniform domains, ambient doubling measures restrict to doubling measures and give zero mass to the boundary; see Murugan~\cite[Lemma~3.5]{Mur}. Zero-angle cusps violate the underlying interior corkscrew condition, yet area remains doubling on a quadratic cusp, as noted by Garc\'ia-Bravo, Ikonen, and Zhu~\cite[Example~5.4]{GIZ}. They also discuss increasing cusp widths for which restricted area is not doubling. Dai, Wen, Xing, and Zhan~\cite{DWXZ} constructed Jordan domains that carry doubling measures but admit no doubling restriction of any ambient doubling measure.

Related work concerns null sets and thickness for doubling measures. Wu~\cite{Wu} studied null sets for doubling and dyadic doubling measures; Buckley, Hanson, and MacManus~\cite{BHM} considered doubling conditions defined by general sets and quasisymmetrically thick Cantor sets. Thin and fat sets were studied by Cs\"ornyei and Suomala~\cite{CS}, Ojala, Rajala, and Suomala~\cite{ORS}, and W.~Wang, S.~Wen, and Z.-Y.~Wen~\cite{WWWfat}. Our domains have positive mass under every full-support ambient doubling measure; the issue is how these restrictions compare with internal doubling measures.

Kaufman and Wu~\cite[p.~537]{KW} also asked which doubling measures on perfect sets arise as ambient restrictions. Complementary mass constructions on Whitney modification sets were studied by Wang, Wen, and Wen~\cite{WWW1,WWW2}; related restriction questions for Cantor measures were treated by Wang and Wen~\cite{WW} and Wei, Wen, and Wen~\cite{WeiWW}. Here the two quantifiers are encoded by the two inclusions between $\mathscr R(\Omega)$ and $\D(\Omega)$. The potential criterion retains the complete prescribed internal measure, rather than only its mass up to comparison.

The low-dimensional compact sets with positive doubling mass used in Theorem~\ref{thm:partialc1} are discussed by Rossi and Shmerkin~\cite[Section~4]{RS}; we supply the needed endpoint-normalized construction. Our finite-system criterion uses consistent masses on nested squares, in the spirit of nested decompositions~\cite{KRS,HK}, and compares all adjacent squares at each level.

Section~2 develops the geometric estimates, projection lemmas, Whitney mass replacement, and a positive-kernel construction. Section~3 proves the stated set relations, the membership criteria, and the examples. The last geometric remark treats finitely many cusps.

\section{Main lemmas}

We use $A\lesssim B$ for $A\le CB$, with constants independent of the points and scales in an estimate; $A\asymp B$ means two-sided comparability. Constants may depend on a fixed domain and on the doubling constants of fixed measures.

We will also use a simple localization observation. A finite cover of $E$ by cusp neighborhoods, ordinary boundary charts, and interior neighborhoods turns uniform local doubling estimates into a doubling estimate at all sufficiently small scales. At the remaining radii $r\ge r_*>0$, a finite cover by relative balls of radius $r_*/4$ gives a uniform positive lower bound for ball masses, while the total mass is finite. Thus these local estimates imply global doubling.

\subsection{Geometric comparisons and ambient estimates}

\begin{lemma}\label{lem:compare}
If $\mu\in\D(X)$ and $\lambda\ge1$, then
\begin{equation}\label{eq:iterate}
\mu(B_X(x,\lambda r))\le D_\mu^{\lceil\log_2\lambda\rceil}\mu(B_X(x,r)).
\end{equation}
For $0<r\le R$, with $s=\log_2D_\mu$, one also has
\begin{equation}\label{eq:lower}
\mu(B_X(x,r))\ge D_\mu^{-1}(r/R)^s\mu(B_X(x,R)).
\end{equation}
Balls of comparable radii with centers separated by a fixed multiple of those radii have comparable masses. Doubling is preserved by bi-Lipschitz maps and by replacing a measure by one comparable to it on all Borel sets.
\end{lemma}
\begin{proof}
Iterate the doubling inequality and use ball inclusions. For a bi-Lipschitz map $T$ with constant $L$, use
$B(x,r/L)\subset T^{-1}(B(Tx,r))\subset B(x,Lr)$.
If $c\mu\le\lambda\le C\mu$, then
$\lambda(B_X(x,2r))\le(C/c)D_\mu\lambda(B_X(x,r))$.
\end{proof}

\begin{lemma}\label{lem:straight}
The local representation \eqref{eq:width}--\eqref{eq:model} holds at the cusp $0$ of a single-cusp piecewise $C^1$ Jordan domain. The straightening map can be chosen as a global bi-Lipschitz homeomorphism of $\R^2$ that preserves area. The relations between $\mathscr R$ and $\D$ are invariant under this change of coordinates.
\end{lemma}
\begin{proof}
The horizontal component of each incident tangent vector is nonzero. The one-sided inverse function theorem therefore gives the two graph representations up to $0$. Their midpoint $h$ is Lipschitz on $[0,a]$. Continue it by $0$ for $t<0$ and by $h(a)+h'(a)(t-a)$ for $t>a$, giving a global $C^1$ Lipschitz function. Define
$T(t,y)=(t,y-h(t))$. Both $T$ and $T^{-1}$ are Lipschitz. Fubini's theorem gives area preservation. For every Borel measure $\nu$,
\[
T_\#(\nu\restr\Omega)=(T_\#\nu)\restr T(\Omega).
\]
Lemma~\ref{lem:compare} proves the assertion about the two measure classes. If the original arcs are smooth up to $0$, the graph functions are smooth there as well.
\end{proof}

\begin{lemma}\label{lem:cork}
Suppose that, for all $x\in E$ and $0<r<r_*$, there is an open ball
\begin{equation}\label{eq:cork}
B(a,\kappa r)\subset\Omega\cap B(x,r)
\end{equation}
with a fixed $\kappa>0$. Then every $\nu\in\D(\R^2)$ has a doubling restriction to $E$. The condition holds locally uniformly at regular $C^1$ boundary points and at junctions of regular $C^1$ arcs with nonzero interior angle, including angle $2\pi$.
\end{lemma}
\begin{proof}
The inclusions $B(x,2r)\subset B(a,3r)$ and \eqref{eq:cork} give
\[
\nu(B(x,2r)\cap E)\le\nu(B(a,3r))
\le D_\nu^{\lceil\log_2(3/\kappa)\rceil}
\nu(B(a,\kappa r))\le C\nu(B(x,r)\cap E).
\]
Larger radii are handled by the localization observation above. A Lipschitz graph chart supplies the interior balls at a regular point. If the tangent rays form an interior angle strictly between $0$ and $2\pi$, a smaller sector between the limiting rays contains such balls at the junction. The $C^1$ approximation to the rays gives the same assertion at nearby centers. When the interior angle is $2\pi$, the excluded region lies between two arcs with a common tangent ray. Balls on the other side of either arc and behind the tip give interior balls of radius comparable to the scale. Equivalently, after scaling around the tip, any fixed sector avoiding that ray lies inside the domain at sufficiently small scales; balls centered away from the tip are treated in a graph chart. A finite cover gives uniform constants on compact boundary portions with no zero-angle junction.
\end{proof}

\begin{lemma}\label{lem:slab}
For $\nu\in\D(\R^2)$ there exist $C\ge1$ and $\eta>0$, depending only on $D_\nu$, such that
\begin{equation}\label{eq:slab}
\nu\big([-r,r]^2\cap\{|y|\le b\}\big)
\le C(b/r)^\eta\nu([-r,r]^2),\qquad 0<b\le r.
\end{equation}
The same estimate holds after translation and rotation, with the same constants.
\end{lemma}
\begin{proof}
Lines are $\nu$-null by the porosity and differentiation argument following Definition~\ref{def:domain}. Divide a square into nine equal subsquares. Comparing an inscribed ball of a subsquare with a circumscribed ball of the parent, doubling gives $\theta=\theta(D_\nu)>0$ such that each subsquare has at least $\theta$ times the parent mass. Retaining only the middle horizontal row discards at least that fraction. After $k$ repetitions the retained set is
\[
[-r,r]\times[-3^{-k}r,3^{-k}r]
\]
and its mass is at most $(1-\theta)^k\nu([-r,r]^2)$. Boundaries contribute zero mass. Take $k=\lfloor\log_3(r/b)\rfloor$ and $\eta=-\log(1-\theta)/\log3$. Rotations and translations preserve the doubling constant.
\end{proof}

\begin{lemma}\label{lem:localobstruction}
Suppose $E\subset\R^2$ is closed and $\mu=\nu\restr E\in\D(E)$ for some $\nu\in\D(\R^2)$. Fix $A\ge1$. If $p\in E$, $0<b,\rho\le r$, and
\begin{equation}\label{eq:localgeometry}
E\cap B(p,r)\subset\{x:\dist(x,L)\le b\},\qquad
|a-p|\le Ar,\qquad B(a,\rho)\subset E,
\end{equation}
where $L$ is a line through $p$, then
\begin{equation}\label{eq:localnecessary}
\frac b r\ge c\left(\frac\rho r\right)^N
\end{equation}
for $c,N>0$ depending only on $A,D_\nu,D_\mu$.
\end{lemma}
\begin{proof}
Let $Q$ be the square of side $2r$ centered at $p$ with one side parallel to $L$. Lemma~\ref{lem:slab} gives
\[
\mu(B_E(p,r))\le C(b/r)^\eta\nu(Q).
\]
Put $s=\log_2D_\nu$. Since $Q\subset B(a,(A+\sqrt2+1)r)$,
\[
\nu(Q)\le C_A(r/\rho)^s\nu(B(a,\rho))
=C_A(r/\rho)^s\mu(B(a,\rho)).
\]
Also $B(a,\rho)\subset B_E(p,(A+1)r)$, so the last mass is at most $C_{A,D_\mu}\mu(B_E(p,r))$. Division by the positive mass of $B_E(p,r)$ yields
$1\le C(b/r)^\eta(r/\rho)^s$, which proves the result with $N=s/\eta$.
\end{proof}

\begin{corollary}\label{cor:local}
If a sequence of configurations \eqref{eq:localgeometry}, with fixed $A$, satisfies
\begin{equation}\label{eq:logobstruction}
\frac{\log(r_k/b_k)}{1+\log(r_k/\rho_k)}\longrightarrow\infty,
\end{equation}
then no full-support ambient doubling measure restricts to a doubling measure on $E$.
\end{corollary}
\begin{proof}
The logarithm of \eqref{eq:localnecessary} bounds the numerator by
$\log(1/c)+N\log(r_k/\rho_k)$, contradicting \eqref{eq:logobstruction}.
\end{proof}

\begin{lemma}\label{lem:graph}
For every Lipschitz function $g:\R\to\R$, the measure
\begin{equation}\label{eq:graphweight}
d\nu_g(t,y)=|y-g(t)|^2\,dt\,dy
\end{equation}
is a full-support doubling measure on $\R^2$.
\end{lemma}
\begin{proof}
In the maximum metric, $d\sigma=dt\,v^2dv$ satisfies
\[
\sigma\big((t-r,t+r)\times(v-r,v+r)\big)
=4v^2r^2+\tfrac43r^4.
\]
Its doubling constant is at most $16$. The maximum and Euclidean metrics are bi-Lipschitz equivalent. The global bi-Lipschitz shear $(t,v)\mapsto(t,v+g(t))$ sends $\sigma$ to $\nu_g$, as follows from Fubini's theorem. Apply Lemma~\ref{lem:compare}. Every nonempty open ball has positive mass.
\end{proof}

\begin{lemma}\label{lem:powerarea}
For $\alpha>1$, restricted area on a domain with local cusp width $t^\alpha$ is doubling near the cusp, uniformly over centers and sufficiently small radii.
\end{lemma}
\begin{proof}
Write $S(p,r)=p+(-r,r)^2$, where $p=(t,y)$ and $|y|\le t^\alpha$. Whenever $t+r$ is sufficiently small,
\begin{equation}\label{eq:powerarea}
\Leb^2(S(p,r)\cap E)\asymp
r\min\{r,(t+r)^\alpha\}.
\end{equation}
For the upper bound, the horizontal length is at most $2r$ and every vertical intersection has length at most $2\min\{r,(t+r)^\alpha\}$. For the lower bound, first suppose $r\le t/2$. For $s\in[t,t+r/2]$, one has $s^\alpha\ge t^\alpha\ge|y|$ and $s^\alpha\asymp(t+r)^\alpha$. The intersection of $(y-r,y+r)$ with $[-s^\alpha,s^\alpha]$ has length at least $\min\{r,s^\alpha\}$. Integration proves the estimate.

If $r>t/2$, then $t+r<3r$. Shrinking the cusp neighborhood ensures $(t+r)^\alpha<r/4$. Consequently $|y|<r/4$, and for $s\in[t+r/4,t+r/2]$ the full vertical section lies inside $(y-r,y+r)$. Its integrated area is at least $cr(t+r)^\alpha$. This proves \eqref{eq:powerarea}, including the case $t=0$.

Since $t+2r\le2(t+r)$, the right side of \eqref{eq:powerarea} changes by a bounded factor when $r$ doubles. Metric equivalence gives Euclidean doubling.
\end{proof}

\begin{lemma}\label{lem:transverse}
Let $f,g\in C^1([0,a])$ vanish at $0$, be positive on $(0,a]$, and satisfy $c\le f/g\le C$. A global bi-Lipschitz homeomorphism of $\R^2$ maps the closed truncated cusp of width $g$ onto that of width $f$, with area distortion bounded above and below. It preserves the relations between $\mathscr R$ and $\D$ and the local doubling of restricted area.
\end{lemma}
\begin{proof}
On the vertical strip $0<t<a$, set $T(t,y)=(t,H(t,y))$, where
\begin{equation}\label{eq:transverse}
H(t,y)=\begin{cases}
yf(t)/g(t),& |y|\le g(t),\\
y+\operatorname{sgn}(y)(f(t)-g(t)),& |y|>g(t).
\end{cases}
\end{equation}
For fixed $t$ this is an increasing piecewise linear homeomorphism, with slopes bounded above and below by positive constants. In the inner region,
\[
|\partial_tH|\le |f'|+C|g'|,
\]
and outside it $|\partial_tH|\le |f'|+|g'|$. These bounds also control variations in $t$ across the matching curves, since the formulas agree there. For example, for fixed $y\ne0$, write
$H=y+\operatorname{sgn}(y)(f-g)\min\{1,|y|/g\}$; on every compact subinterval of $(0,a)$ this is Lipschitz and satisfies the stated derivative bound almost everywhere. Integrating, then passing to $t=0$, gives a uniform bound on the full interval. Set $H(0,y)=y$. Interchanging $f$ and $g$ gives the same estimates for the inverse. Thus $T$ is bi-Lipschitz up to the tip.

One may continue $T$ by the identity for $t<0$ and by its fixed vertical formula at $a$ for $t>a$, obtaining a global bi-Lipschitz map. Its vertical Jacobian is $f/g$ inside the cusp and $1$ outside. Fubini's theorem, or the area formula, shows that the pushforward of area is comparable to area on all Borel sets. The identity for restrictions in Lemma~\ref{lem:straight} and Lemma~\ref{lem:compare} complete the proof.
\end{proof}

\subsection{Axial projections and reflected masses}

\begin{lemma}\label{lem:projection}
Let $f$ be Lipschitz on $[0,1]$, positive on $(0,1]$, and satisfy $f(0)=0$. If $\mu\in\D(E_f)$, then the axial projection $\lambda=\pi_\#\mu$, $\pi(t,y)=t$, is a full-support doubling measure on $[0,1]$, with constant controlled by $D_\mu$ and $\Lip(f)$.
\end{lemma}
\begin{proof}
Write $I(t,r)=(t-r,t+r)\cap[0,1]$. If $r\ge cf(t)$ for fixed $c>0$, then
$\pi^{-1}I(t,2r)\subset B_{E_f}((t,0),C_cr)$, since $f(s)\le f(t)+2\Lip(f)r$. Also $B_{E_f}((t,0),r)\subset\pi^{-1}I(t,r)$, and doubling gives the desired inequality.

If $r<cf(t)$, choose points $p_j$ along the cross-section at $t$ with spacing comparable to $r$. Lipschitz regularity makes $B_{E_f}(p_j,Cr)$ cover $\pi^{-1}I(t,2r)$. For fixed small $b>0$, the balls $B_{E_f}(p_j,br)$ are disjoint and lie in $\pi^{-1}I(t,r)$. Apply doubling to each pair and sum:
\[
 \lambda(I(t,2r))\le\sum_j\mu(B_{E_f}(p_j,Cr))
 \le C\sum_j\mu(B_{E_f}(p_j,br))\le C\lambda(I(t,r)).
\]
At the endpoints, truncate the intervals; the argument is unchanged.
\end{proof}
\begin{lemma}\label{lem:canonical}
The measure $\eta_\lambda$ in \eqref{eq:canonicalmeasure} is doubling on $E_f$, and its axial projection is $\lambda$. Near the cusp its ball masses satisfy \eqref{eq:canonicalball}.
\end{lemma}
\begin{proof}
A full-support doubling measure on a nondegenerate interval has no atoms, by the separated-ball estimate in the boundary argument following Definition~\ref{def:domain}, so the undefined value of the density in \eqref{eq:canonicalmeasure} at $t=0$ is immaterial. Integration over the cross-sections gives $\pi_\#\eta_\lambda=\lambda$.

We first prove that the canonical measure is doubling. Shrink the cusp neighborhood so that $|f'|$ is sufficiently small, and write $p=(t,y)\in E_f$. If $r\le c f(t)$, then $f(s)\asymp f(t)$ for $|s-t|\le r$. Inscribed and circumscribed rectangles, together with $|y|\le f(t)$, give
\[
 \eta_\lambda(B_{E_f}(p,r))
 \asymp \frac r{f(t)}\lambda(I(t,r)).
\]
If $r\ge C f(t)$, an inscribed rectangle contains full cross-sections over an axial interval of length comparable to $r$. The ball mass is therefore comparable to $\lambda(I(t,r))$. The intermediate range is controlled by the same rectangle estimates. Thus, uniformly near the cusp,
\begin{equation}\label{eq:canonicalball}
 \eta_\lambda(B_{E_f}(p,r))
 \asymp\min\{1,r/f(t)\}\lambda(I(t,r)),
\end{equation}
where the minimum is interpreted as $1$ when $t=0$. Doubling $r$ changes the right-hand side by at most a fixed factor. Away from the cusp, $f$ has a positive lower bound, and the local density is comparable to that of the product measure $d\lambda(t)dy$. The rectangle estimates for boundary graphs apply there as well. Compactness handles the remaining positive scales, proving that $\eta_\lambda$ is doubling.

\end{proof}
\begin{lemma}\label{lem:corkrealization}
Suppose $E=\overline\Omega$ is bounded and satisfies a uniform interior ball condition: for some $\kappa,r_*>0$, every $x\in\partial E$ and $0<r<r_*$ admit
\[
B(a,\kappa r)\subset\Omega\cap B(x,r).
\]
Then $\D(E)\subseteq\mathscr R(E)$. For each $\mu\in\D(E)$, the representing measure $\nu\in\D(\R^2)$ can be chosen absolutely continuous outside $E$.
\end{lemma}
\begin{proof}
Use a Whitney decomposition of $\R^2\setminus E$, with side lengths $s_Q$ comparable to distance from $E$. For every sufficiently small square, apply the interior ball condition at a nearest boundary point and at a fixed multiple of $s_Q$. Shrinking the resulting ball by a factor of two gives $A_Q\subset\Omega$ such that its radius, its boundary depth, and its distance from $Q$ are comparable to $s_Q$. Thus the balls $A_Q$ have bounded overlap across all scales: a point in $A_Q$ has distance from $\partial E$ comparable to $s_Q$, and only boundedly many squares of a fixed scale can occur nearby.

Assign $m_Q=\mu(A_Q)$ to these squares and $m_Q=s_Q^2$ to all sufficiently distant squares, with arbitrary positive masses on the finitely many intermediate squares. Set
\[
\nu=\mu+\sum_Q\frac{m_Q}{s_Q^2}\Leb^2\restr Q.
\]
The overlap estimate gives local finiteness. Adjacent reflected squares have comparable masses by internal doubling, while the other interfaces involve only finitely many positive masses. Hence adjacent densities are uniformly comparable.

For $x\in E$ and small $r$, every square meeting $B(x,r)$ has $s_Q\lesssim r$ and $A_Q\subset B_E(x,Cr)$. Consequently,
\[
\mu(B_E(x,r))\le\nu(B(x,r))
\le \mu(B_E(x,r))+C\mu(B_E(x,Cr))
\lesssim\mu(B_E(x,r)).
\]
This proves doubling for centers in $E$. For $z\notin E$, put $d=\dist(z,E)$ and choose a nearest $q\in E$. If $d\ge4r$, local Whitney densities are comparable throughout $B(z,2r)$. If $d<r/2$, compare with balls centered at $q$. If $r/2\le d<4r$, squares meeting $B(z,r/8)$ have side comparable to $r$ and assigned mass comparable to $\mu(B_E(q,r))$. Integration on this smaller ball gives the lower bound for $\nu(B(z,r))$, while $B(z,2r)\subset B(q,6r)$ gives the upper bound.

Uniform small-scale doubling follows. On a compact set of centers and an interval of radii bounded away from zero, local finiteness and a finite ball cover give uniform upper and positive lower bounds. The density equals $1$ outside a fixed large ball, so sufficiently large balls have mass comparable to their area. This proves global doubling, and the construction gives $\nu\restr E=\mu$.
\end{proof}
\subsection{A Whitney replacement principle}
\begin{lemma}\label{lem:whitneyreplacement}
Let $f\in C^1([0,1])$ satisfy $f(0)=f'(0)=0$ and $f(t)>0$ for $t>0$,
and put $E=E_f$. Let $\mu\in\D(E)$ and $\lambda=\pi_\#\mu$.
Suppose that $\sigma\in\D(\R^2)$ has the following properties in a sufficiently small cusp neighborhood:
\begin{align}
 \sigma(B((t,0),f(t)))&\asymp\lambda(I(t,f(t))),\label{eq:replacewidth}\\
 \sigma(B(p,r))&\gtrsim\lambda(I(t,r)),
 \qquad p=(t,y)\in E,\quad r\ge f(t),\quad t+r<r_*.
 \label{eq:replacedominance}
\end{align}
At $t=0$ only the second condition is required.
Then $\mu\in\mathscr R(E)$, with a representing measure $\nu\in\D(\R^2)$ that is absolutely continuous on $\R^2\setminus E$.
\end{lemma}
\begin{proof}
We give the construction and the comparisons needed at its interfaces.
The projection estimates imply, for every fixed $c>0$,
\begin{equation}\label{eq:replacestrip}
 \mu(B_E(p,r))\asymp\lambda(I(t,r))
 \qquad(p=(t,y)\in E,\ r\ge cf(t),\ t+r<r_*).
\end{equation}
Indeed, $B_E(p,r)$ lies in the full-section strip over $I(t,r)$, while that strip is contained in $B_E(p,C_cr)$ because $|y|\le f(t)$ and
$f(s)\le f(t)+\|f'\|_\infty|s-t|$.
Doubling of $\mu$ gives the comparison.
Moving centers and changing radii by fixed factors in \eqref{eq:replacewidth} therefore gives
\begin{equation}\label{eq:replaceinterface}
 \sigma(B(p,r))\asymp\mu(B_E(p,r))\asymp\lambda(I(t,r))
 \quad\text{if }c_1f(t)\le r\le c_2f(t),
\end{equation}
for any fixed $0<c_1<c_2<\infty$.
In particular, \eqref{eq:replacedominance} remains valid with $r\ge cf(t)$ for any fixed $c>0$.

\emph{The mass assignment.}
Take a Whitney decomposition $\mathcal W$ of $\R^2\setminus E$ into squares with disjoint interiors.
Write $s_Q$ for the side of $Q$ and $(u_Q,v_Q)$ for its center. We use the usual properties
\[
 s_Q\asymp\dist(Q,E),\qquad
 s_Q\asymp s_{Q'}\quad\text{if }\overline Q\cap\overline {Q'}\ne\varnothing.
\]
All constants implicit in the Whitney properties are fixed first.
Choose a small number $b>0$, then shrink the cusp chart until its slope bound $L$ is sufficiently small in terms of these constants and $b$.
For squares in this chart set $T_Q=\max\{u_Q,s_Q\}$.
Call a square fine if $s_Q<bf(T_Q)$, and coarse otherwise.
For fine squares $T_Q=u_Q$ and $u_Q\gg s_Q$: the alternative $T_Q=s_Q$ is excluded by $f(s)=o(s)$.
If $(\tau,\pm f(\tau))$ is a nearest boundary point, then $|\tau-u_Q|\lesssim s_Q$ and
\[
 f(\tau)\ge f(u_Q)-CLs_Q\gg s_Q.
\]
An inward vertical displacement of a fixed multiple of $s_Q$ thus supplies an interior ball $A_Q$ with
\begin{equation}\label{eq:reflectdepth}
 \operatorname{rad}(A_Q)\asymp s_Q,\qquad
 \dist(A_Q,\partial E)\asymp s_Q,\qquad
 \dist(A_Q,Q)\lesssim s_Q.
\end{equation}
The same construction is available for sufficiently small squares next to the remaining compact boundary, using its finite graph and nonzero-angle charts.

These reflected balls have bounded overlap, including across different scales.
If $z\in A_Q$, then \eqref{eq:reflectdepth} gives $s_Q\asymp\dist(z,\partial E)$ and $|z-(u_Q,v_Q)|\lesssim s_Q$.
Only a fixed number of dyadic scales are possible, and at each such scale disjointness of Whitney interiors bounds the number of possible squares.
Consequently $\sum_Q\mathbf1_{A_Q}\le C$ for all reflected balls.

Assign $m_Q=\mu(A_Q)$ to the fine squares and $m_Q=\sigma(Q)$ to the coarse squares in the cusp chart.
Away from the cusp assign reflected-ball masses to all sufficiently small boundary squares.
Set $m_Q=as_Q^2$, for one fixed $a>0$, on all squares sufficiently far from $E$; the finitely many remaining intermediate squares may be given arbitrary positive masses.
This prescription can be made consistently on overlapping charts: along a compact interface away from $0$, $f$ has a positive lower bound, so all sufficiently small squares are reflected on both sides. Any interface squares to which different rules apply have a fixed positive minimum size, and there are only finitely many of them.
Define
\begin{equation}\label{eq:replacementformula}
 \nu=\mu+\sum_{Q\in\mathcal W}\frac{m_Q}{s_Q^2}\Leb^2\restr Q.
\end{equation}
Boundary conventions for the squares do not matter for the area terms.
Bounded overlap controls the sum of reflected masses by $C\mu(E)$.
The coarse masses in each compact set are controlled by the $\sigma$-mass of a fixed enlargement of that set.
The remaining terms are locally finite. Hence $\nu$ is a full-support Radon measure and $\nu\restr E=\mu$ exactly.

\emph{Adjacent squares.}
Two adjacent fine squares have comparable assigned masses by doubling of $\mu$ and \eqref{eq:reflectdepth}; two adjacent coarse squares have comparable masses by doubling of $\sigma$.
At a mixed pair, comparability of the side lengths and
$|f(T_Q)-f(T_{Q'})|\le L|T_Q-T_{Q'}|\lesssim Ls_Q$
give $s_Q\asymp f(T_Q)$, with constants allowed to depend on $b$.
Equation~\eqref{eq:replaceinterface}, the strip comparison, and doubling at nearby centers give
\[
 \mu(A_Q)\asymp\lambda(I(T_Q,s_Q))\asymp\sigma(Q').
\]
Thus adjacent exterior densities $m_Q/s_Q^2$ are uniformly comparable.
The ordinary reflected rules give the same comparison elsewhere; the finite interfaces only enlarge the constant.

\emph{Centers in $E$.}
Fix $p=(t,y)\in E$ and work at scales with $t+r$ smaller than a fixed subchart radius.
If $r\le cf(t)$ for a sufficiently small fixed $c>0$, every exterior square meeting $B(p,r)$ is fine, and its reflected ball lies in $B_E(p,Cr)$.
Bounded overlap and doubling give
\begin{equation}\label{eq:replacementfine}
 \mu(B_E(p,r))\le\nu(B(p,r))\le C\mu(B_E(p,r)).
\end{equation}
The first inequality also covers balls lying entirely inside $E$.

For $r\ge cf(t)$, the interior mass and the sum of the fine-square masses are bounded by
$C\mu(B_E(p,Cr))\asymp\lambda(I(t,r))$.
Every square meeting $B(p,r)$ has side at most $Cr$ and lies in $B(p,Cr)$.
The coarse-square sum is consequently at most $\sigma(B(p,Cr))$.
Using \eqref{eq:replacedominance} gives
$\nu(B(p,r))\lesssim\sigma(B(p,r))$.
If $r/f(t)$ stays in a fixed bounded interval, the reverse inequality follows from
$\nu\ge\mu$ and \eqref{eq:replaceinterface}.
Otherwise, choose a ball $U\subset B(p,r)$ of radius $r/16$ whose height above the axis is approximately $r/2$.
Since $f(t)\ll r$ and the slope is small, $\dist(U,E)\asymp r$.
Every Whitney square meeting $U$ is coarse, with side comparable to $r$.
For each such square, doubling of $\sigma$ gives
$\sigma(Q)/s_Q^2\asymp\sigma(B(p,r))/r^2$.
Integrating over $U$ proves
\begin{equation}\label{eq:replacementcoarseprofile}
 \nu(B(p,r))\asymp\sigma(B(p,r)),\qquad r\ge cf(t).
\end{equation}
This argument includes $p=0$, because $f(s)=o(s)$.
Equations~\eqref{eq:replacementfine}--\eqref{eq:replacementcoarseprofile}, with
\eqref{eq:replaceinterface} in their fixed-ratio transition, prove doubling for centers in $E$ near the cusp.
Away from it, the reflected construction gives \eqref{eq:replacementfine} at all sufficiently small scales, by the same bounded-overlap estimate.

\emph{Exterior centers.}
Let $z\notin E$, $d=\dist(z,E)$, and choose a nearest point $q\in E$.
If $d\ge4r$, the Whitney squares meeting $B(z,2r)$ have comparable sizes; their densities are comparable by a bounded chain of adjacent squares. The desired doubling estimate follows by comparison with area.
If $d<r/2$, use
\[
 B(q,r/2)\subset B(z,r)\subset B(z,2r)\subset B(q,3r)
\]
and doubling at the center $q$.

Suppose $r/2\le d<4r$.
Every square meeting $B(z,r/8)$ has side comparable to $r$ and lies within distance $Cr$ of $q$.
For a fine square, \eqref{eq:reflectdepth} gives
$m_Q\asymp\mu(B_E(q,r))$; its fine-scale condition and the small-slope estimate imply that $r/f(t_q)$ is bounded above by a fixed constant.
Equations~\eqref{eq:replacementfine}, \eqref{eq:replacementcoarseprofile}, and
\eqref{eq:replaceinterface} therefore identify this mass with $\nu(B(q,r))$.
For a coarse square, $r\gtrsim f(t_q)$, with a fixed constant allowed to depend on $b$, so
$m_Q=\sigma(Q)\asymp\sigma(B(q,r))\asymp\nu(B(q,r))$.
At $t_q=0$ only this latter case is possible at sufficiently small scales.
It follows that the density throughout $B(z,r/8)$ is comparable to
$\nu(B(q,r))/r^2$. Thus
\[
 \nu(B(z,r))\gtrsim\nu(B(q,r)),\qquad
 \nu(B(z,2r))\le\nu(B(q,6r))\lesssim\nu(B(q,r)).
\]
The same reasoning with reflected balls handles exterior centers near the ordinary boundary.

We have established uniform doubling at all sufficiently small scales in a fixed neighborhood of $E$.
On any remaining compact set of centers and compact interval of positive radii, local finiteness supplies a uniform upper bound and a finite ball cover together with full support supplies a positive lower bound.
Finally, $\nu$ has the constant density $a$ outside a fixed large ball.
Balls contained there satisfy doubling directly; for the other sufficiently large balls, the finite mass of the exceptional compact set and an exterior subball of comparable radius give $\nu(B(z,r))\asymp r^2$.
These observations prove global doubling and complete the construction.
\end{proof}
\subsection{A positive-kernel lift}
\begin{lemma}\label{lem:poissonlift}
Let $\rho$ be a nonzero finite positive Borel measure on $\mathbb R$ with compact support, and let $0<\varepsilon<1$ and $N\ge2$.
Define
\[
 M_\rho(x,h)=\int\left(\frac h{h+|x-s|}\right)^N d\rho(s),
 \qquad
 d\sigma(x,y)=|y|^{\varepsilon-2}M_\rho(x,|y|)\,dx\,dy
 \quad(y\ne0),
\]
and assign zero mass to the axis. Then $\sigma$ is a doubling Radon measure with full support on $\mathbb R^2$, and
\begin{equation}\label{eq:liftprofile}
 \sigma\big((t-r,t+r)\times(-r,r)\big)
 \asymp r^\varepsilon M_\rho(t,r).
\end{equation}
The comparison constants depend only on $N$ and $\varepsilon$.
\end{lemma}
\begin{proof}
By Tonelli's theorem, the mass of a square centered on the axis equals
\[
 2\int d\rho(s)\int_0^r h^{\varepsilon-2}
   \int_{t-r}^{t+r}\left(\frac h{h+|x-s|}\right)^N dx\,dh.
\]
If $|s-t|\le2r$, the inner $x$-integral is at most $C_Nh$, so this contribution is bounded by
$C r^\varepsilon\rho([t-2r,t+2r])$.
If $R=|s-t|>2r$, the inner integral is at most $Cr(h/R)^N$. Its integral in $h$ is at most
$C r^\varepsilon(r/R)^N$, since $N+\varepsilon>1$.
Adding the two contributions gives the upper bound in \eqref{eq:liftprofile}.
For the lower bound, restrict to $h\in[r/2,r]$. Whenever $|x-t|\le r$, we then have
$M_\rho(x,h)\asymp M_\rho(t,r)$, which proves the lower estimate.
This also proves local finiteness, including at points where $\rho$ may have atoms on the axis.

Pointwise comparison of the kernels gives
$M_\rho(t,2r)\le2^N M_\rho(t,r)$, so balls centered on the axis satisfy a uniform doubling estimate.
On balls whose radius is at most one quarter of the distance of their center from the axis, the density is uniformly comparable to its value at the center, and doubling follows directly.
For the remaining balls, use a containing ball centered on the axis for the upper bound, and an interior ball away from the axis, with comparable radius, for the lower bound.
The mass of the latter is comparable to the same quantity $r^\varepsilon M_\rho(t,r)$, which proves doubling in this case as well.
The kernel is strictly positive, so $\sigma$ has full support.
\end{proof}

\section{Proofs of the theorems}

\subsection{A bad restriction at every zero-angle cusp}

\begin{proof}[Proof of Theorem~\ref{thm:bad}]
Use Lemma~\ref{lem:straight} to straighten the cusp, retaining the notation $\Omega$ for the resulting domain. Put $m=\Leb^2\restr E$. If $m$ is not doubling, ambient area itself gives the conclusion. Otherwise let $D$ be its doubling constant.

Take $t_k=8^{-k}$ and $r_k=t_k/8$, starting sufficiently far into the cusp chart. Choose a fixed piecewise linear function $\psi:\R\to[0,1]$, supported on $[5/4,7/4]$, equal to $1$ on $[11/8,13/8]$, and linear on each transition interval. Define
\begin{equation}\label{eq:g}
g(t)=\sum_{k\ge k_0}\frac{r_k}{8}\,
\psi\left(\frac{t-t_k}{r_k}\right).
\end{equation}
The supports are disjoint and the slopes are uniformly bounded. Their heights tend to zero at the only accumulation point, so $g$ is globally Lipschitz, with $g(0)=0$. Set
\[
p_k=(t_k,0),\quad q_k=(t_k+3r_k/2,0),\quad
a_k=\sup_{t\in[t_k-r_k,t_k+2r_k]}f(t).
\]
The condition $f(t)=o(t)$ gives $a_k/r_k\to0$. Increase $k_0$ so that $a_k\le r_k/16$.

All perturbations vanish on $B(p_k,r_k)$. For the ambient measure \eqref{eq:graphweight},
\begin{equation}\label{eq:smallmass}
\nu_g(B(p_k,r_k)\cap E)\le a_k^2m(B_E(p_k,r_k)).
\end{equation}
The ball $B(q_k,r_k/16)$ lies in $B(p_k,2r_k)$, and its horizontal coordinates lie in the plateau where $g=r_k/8$. On its intersection with $E$ we have $|y-g(t)|\ge r_k/16$. Also
$B_E(p_k,r_k)\subset B_E(q_k,3r_k)$. Since $48<2^6$,
\begin{equation}\label{eq:largemass}
\begin{split}
\nu_g(B(p_k,2r_k)\cap E)
&\ge\frac{r_k^2}{256}m(B_E(q_k,r_k/16))\\
&\ge\frac{r_k^2}{256D^6}m(B_E(p_k,r_k)).
\end{split}
\end{equation}
Combining \eqref{eq:smallmass} and \eqref{eq:largemass} gives
\begin{equation}\label{eq:badquotient}
\frac{\nu_g(B(p_k,2r_k)\cap E)}{\nu_g(B(p_k,r_k)\cap E)}
\ge\frac1{256D^6}\left(\frac{r_k}{a_k}\right)^2\longrightarrow\infty.
\end{equation}
The ambient measure is doubling by Lemma~\ref{lem:graph}, but its restriction is not. Pulling it back under the original shear preserves doubling and absolute continuity. Finally apply \eqref{eq:openclosed}.
\end{proof}

\begin{remark}
When $f(t)\asymp t^\alpha$, the quotients in \eqref{eq:badquotient} grow at least as $c t_k^{-2(\alpha-1)}$.
\end{remark}

\subsection{Infinitely flat cusps}

\begin{proof}[Proof of Theorem~\ref{thm:flat}]
Assume that $\nu\restr\Omega$ is doubling. Straighten the cusp and put $\mu=\nu\restr E$, which is doubling by \eqref{eq:openclosed}. Fix a sufficiently small $R>0$. Lemma~\ref{lem:compare} gives
\begin{equation}\label{eq:tiplower}
\mu(B_E(0,r))\ge c_0r^s\qquad(0<r<R),\quad s=\log_2D_\mu>0.
\end{equation}
Since $f=o(t)$, the envelope $F(r)=\sup_{0\le t\le r}f(t)$ satisfies $F(r)=o(r)$. For sufficiently small $r$,
\[
E\cap B(0,r)\subset[-r,r]^2\cap\{|y|\le F(r)\}.
\]
Lemma~\ref{lem:slab} and local finiteness of $\nu$ imply
\begin{equation}\label{eq:tipupper}
\mu(B_E(0,r))\le C(F(r)/r)^\eta\nu([-r,r]^2)
\le C\nu([-R,R]^2)(F(r)/r)^\eta.
\end{equation}
Together with \eqref{eq:tiplower}, this gives $F(r)\ge cr^{1+s/\eta}$ and proves \eqref{eq:envelopebound}. If $f$ is infinitely flat, then $F(r)\le C_Mr^M$ for each $M$, contradicting this lower bound when $M>1+s/\eta$. If $f$ is nondecreasing, $F(r)=f(r)$.
\end{proof}

\begin{remark}
For nonmonotone widths, isolated small values of $f$ do not control the width of an entire tip neighborhood, so the envelope $F$ is needed. The width $f(t)=e^{-1/t^\beta}$, $\beta>0$, is infinitely flat and gives $\mathscr R(\Omega)\cap\D(\Omega)=\varnothing$.
\end{remark}

\subsection{A finite-logarithmic-order counterexample}

\begin{proof}[Proof of Theorem~\ref{thm:counter}]
Fix $\alpha>1$. Choose a smooth nondecreasing function $\psi:\R\to[0,1]$ which is zero on $(-\infty,0]$ and one on $[1,\infty)$. Put
\begin{equation}\label{eq:stepdata}
t_k=2^{-3^k},\quad L_k=\log(1/t_k),\quad
H_k=t_k^\alpha e^{\sqrt{L_k}},\quad r_k=H_kL_k.
\end{equation}
For $k_0$ sufficiently large, define on a small interval $[0,a]$, with $a<1$,
\begin{equation}\label{eq:stepwidth}
f(t)=t^\alpha+\sum_{k\ge k_0}H_k\psi\left(\frac{t-t_k}{r_k}\right),
\qquad f(0)=0.
\end{equation}
The following elementary limits will be used:
\begin{equation}\label{eq:steplimits}
\frac{r_k}{t_k}\longrightarrow0,\qquad
\frac{H_{k+1}}{H_k}\longrightarrow0,\qquad
\frac{H_{k+1}}{t_k^\alpha}\longrightarrow0,
\qquad \sum_k\frac1{L_k}<\infty.
\end{equation}
Indeed, $L_{k+1}=3L_k$, so the first ratio is
$L_ke^{-(\alpha-1)L_k+\sqrt{L_k}}$; the second and third have exponents
$-2\alpha L_k+(\sqrt3-1)\sqrt{L_k}$ and
$-2\alpha L_k+\sqrt{3L_k}$, respectively. Increase $k_0$ so that
\begin{equation}\label{eq:tail}
r_k<t_k/16,\qquad
\sum_{j>k}H_j\le2H_{k+1},\qquad
\sum_{j\ge k}H_j\le2H_k.
\end{equation}

\emph{Regularity and order.}
The series in \eqref{eq:stepwidth} converges uniformly. Its derivative series converges uniformly as well, because the $k$th derivative term has norm at most $\|\psi'\|_\infty/L_k$. Each derivative term vanishes at zero. Hence $f\in C^1([0,a])$ and $f'(0)=0$. On any compact subinterval of $(0,a]$, all sufficiently small steps are constant and only finitely many derivatives remain. Thus $f$ is smooth at every positive parameter. Its derivative is positive there, since the derivative of $t^\alpha$ is positive and all step derivatives are nonnegative.

The function $G(t)=t^\alpha e^{\sqrt{\log(1/t)}}$ is increasing for sufficiently small $t$. For such $t$, let $j$ be the first index with $t_j\le t$. All terms with smaller index vanish, so \eqref{eq:tail} gives
\[
t^\alpha\le f(t)\le t^\alpha+2H_j
\le t^\alpha+2G(t)\le3G(t).
\]
Choose $a$ small enough that this holds throughout $(0,a)$. This proves \eqref{eq:counterbounds}. With $L=\log(1/t)$, division of the logarithmic bounds by $\log t=-L$ gives
\[
\alpha-L^{-1/2}-\frac{\log3}{L}
\le\frac{\log f(t)}{\log t}\le\alpha.
\]
The logarithmic order is therefore exactly $\alpha$.

\emph{Thin sections.}
Let
\begin{equation}\label{eq:badgeometry}
p_k=(t_k-2r_k,0),\quad a_k=(t_k+2r_k,0),\quad
b_k=2t_k^\alpha,\quad \rho_k=H_k/4.
\end{equation}
For $t\in[t_k-3r_k,t_k-r_k]$, the $k$th step and all earlier steps have not started. Later steps contribute at most $2H_{k+1}=o(t_k^\alpha)$. Since $t<t_k$, for all sufficiently large $k$,
\[
f(t)\le t_k^\alpha+2H_{k+1}\le b_k.
\]
Thus $E\cap B(p_k,r_k)\subset\{|y|\le b_k\}$.

\emph{Nearby interior balls.}
On $B(a_k,\rho_k)$ the horizontal coordinate is larger than $t_k+r_k$, because $\rho_k=H_k/4<r_k$. The $k$th step is complete, so $f(t)\ge H_k>\rho_k$ there. Consequently
\[
B(a_k,\rho_k)\subset\Omega_f,\qquad |a_k-p_k|=4r_k.
\]
For large $k$ one has $0<b_k,\rho_k\le r_k$. All the hypotheses of Lemma~\ref{lem:localobstruction} hold with $A=4$. On the other hand,
\begin{equation}\label{eq:divergence}
\frac{\log(r_k/b_k)}{1+\log(r_k/\rho_k)}
=\frac{\sqrt{L_k}+\log L_k-\log2}{1+\log(4L_k)}
\longrightarrow\infty.
\end{equation}
Corollary~\ref{cor:local} excludes every ambient doubling restriction.

\emph{Completion of the domain.}
The graphs $y=\pm f(t)$ are regular $C^1$ arcs up to $0$. Outside a smaller fixed neighborhood of $0$, smoothly modify the width to a positive constant. Join the resulting horizontal upper and lower ends by a right semicircle. The tangents match at the two joins, so the resulting Jordan boundary is regular $C^1$ away from $0$, and $0$ is its unique zero-angle cusp. All the configurations above lie in the unchanged neighborhood for large $k$. Thus the resulting domain has $\mathscr R(\Omega)\cap\D(\Omega)=\varnothing$. In contrast, a regular cap applied to $f(t)=t^\alpha$ gives $\D(\Omega)\subsetneq\mathscr R(\Omega)$ by Theorem~\ref{thm:power}.
\end{proof}

\begin{remark}
The example is not power-comparable of order $\alpha$. At $t=t_k+2r_k$, one has $t/t_k\to1$ and $f(t)\ge H_k$, so $f(t)/t^\alpha\to\infty$ along these points. For integer $\alpha$, it cannot have smooth endpoint branches either: a smooth positive width of finite logarithmic order $\alpha$ would have leading term $b t^\alpha$ by Taylor's theorem, contrary to this unbounded ratio.
\end{remark}

\subsection{Internal doubling measures that are not ambient restrictions}

\begin{proof}[Proof of Corollary~\ref{thm:internalnonrestriction}]
We work after straightening; the global bi-Lipschitz map $T$ preserves
doubling and the relation of being a restriction. Suppose $\mu=\nu\restr\Omega$, and use the same symbols
for the transported measures. Shrink $a$ so that $f(t)<a$ on $(0,a)$. The thin-strip estimate of
Lemma~\ref{lem:slab}, applied in the fixed square $[-a,a]^2$, gives constants
$C>0$ and $\eta>0$ such that, for $0<h<a$,
\[
\lambda\{t:f(t)<h\}
\le \nu\bigl([-a,a]^2\cap\{|y|<h\}\bigr)
\le C h^\eta.
\]
For $0<\varepsilon<\eta$, Tonelli's theorem yields
\[
\int_{(0,a)}f(t)^{-\varepsilon}\,d\lambda(t)
=a^{-\varepsilon}\lambda((0,a))
 +\varepsilon\int_0^a h^{-\varepsilon-1}
     \lambda\{t:f(t)<h\}\,dh<\infty.
\]
\end{proof}

\begin{remark}\label{rem:momentnonsufficient}
The necessary condition~\eqref{eq:negativewidthmoment} is not sufficient. On each domain in Theorem~\ref{thm:counter}, every $\mu\in\D(\Omega)$ satisfies that condition for some $\varepsilon>0$, yet none is a planar doubling restriction.
\end{remark}
\begin{proof}
Fix a domain from Theorem~\ref{thm:counter} and $\mu\in\D(\Omega)$. Use the same notation for its transported measure on the closure with zero boundary mass. Iterating \eqref{eq:nonatomiccontraction} gives
\[
\mu(B(0,r))\le Cr^\delta
\]
for some $\delta>0$ and all sufficiently small $r$.
Since $f(t)=o(t)$, the local strip with $0<t\le r$ lies in $B(0,2r)$ when
$r$ is small. Hence $\lambda((0,r])\le C r^\delta$. The lower bound
$f(t)\ge t^\alpha$ in~\eqref{eq:counterbounds} now implies, whenever
$0<\varepsilon<\delta/\alpha$,
\[
\int_{(0,a)}f(t)^{-\varepsilon}\,d\lambda(t)
\le\int_{(0,a)}t^{-\alpha\varepsilon}\,d\lambda(t)<\infty.
\]
The last integral is finite by the layer-cake formula. Theorem~\ref{thm:counter} gives disjointness of the two classes.
\end{proof}

\subsection{The potential criterion and axial reduction}
\begin{proof}[Proof of Theorem~\ref{thm:potentialcriterion}: necessity]
Let $\nu\in\mathcal D(\mathbb R^2)$ satisfy $\nu|_{E_f}=\mu$.
Decay in strips around a line gives an exponent $\eta>0$ such that, whenever $0<h\le R$,
\[
 \nu\big((t-R,t+R)\times(-h,h)\big)
 \le C(h/R)^\eta\nu(B((t,0),3R)).
\]
Choose $0<\varepsilon<\min\{1,\eta\}$ and write $d\rho=f^{-\varepsilon}\,d\lambda$.
If $R\ge f(t)$, the Lipschitz property gives $f(s)\le (1+\operatorname{Lip}f)R$ for $s\in I(t,R)$.
The layer-cake formula therefore yields
\begin{align*}
 \rho(I(t,R))
 &\le R^{-\varepsilon}\lambda(I(t,R))
   +\varepsilon\int_0^R h^{-\varepsilon-1}
      \lambda\{s\in I(t,R):f(s)<h\}\,dh\\
 &\le C R^{-\varepsilon}\nu(B((t,0),C_f R)).
\end{align*}
The first term is also bounded by the right-hand side, since the corresponding cross-sections lie in a fixed dilate of the ball.
Choosing one fixed pair $t,R$ whose interval covers the entire axis segment gives $\rho([0,1])<\infty$.

Set $r=f(t)$, with $t$ sufficiently small. Since $f(t)=o(t)$ and the slope is small,
$B((t,0),c r)\subset E_f$, and
\[
 \nu(B((t,0),r))\le C\mu(B((t,0),c r))
 \le C\lambda(I(t,r)).
\]
Writing $d_\nu=\log_2D_\nu$, we obtain, for $R\ge r$,
\[
 \rho(I(t,R))\le C R^{-\varepsilon}(R/r)^{d_\nu}\lambda(I(t,r)).
\]
Decompose the potential integral into $|s-t|<r$ and the annuli $2^kr\le |s-t|<2^{k+1}r$. Choose
$N>\max\{2,d_\nu-\varepsilon\}$. Then
\[
 \int\left(\frac r{r+|t-s|}\right)^N d\rho(s)
 \le C r^{-\varepsilon}\lambda(I(t,r))
       \sum_{k\ge0}2^{-k(N+\varepsilon-d_\nu)}.
\]
This proves \eqref{eq:potentialtesting}.
\end{proof}

\begin{proof}[Proof of Theorem~\ref{thm:potentialcriterion}: sufficiency]
Set $d\rho=f^{-\varepsilon}\,d\lambda$. The test at any one sufficiently small $t>0$ implies $\rho([0,1])<\infty$, because the kernel has a positive lower bound on this compact interval. Let $\sigma$ be the measure in Lemma~\ref{lem:poissonlift}.
Write $r=f(t)$. Near the tip the small-slope estimate gives $f(s)\asymp f(t)$ whenever $|s-t|\le r$. Therefore,
\[
M_\rho(t,r)\gtrsim r^{-\varepsilon}\lambda(I(t,r)).
\]
Together with \eqref{eq:potentialtesting} and \eqref{eq:liftprofile}, this yields
\begin{equation}\label{eq:widthreference}
\sigma(B((t,0),f(t)))\asymp\lambda(I(t,f(t))).
\end{equation}
For $r\ge f(t)$ and $s\in I(t,r)$, the Lipschitz estimate gives $f(s)\le(1+\Lip f)r$. Hence
\[
r^\varepsilon M_\rho(t,r)
\gtrsim r^\varepsilon\int_{I(t,r)}f(s)^{-\varepsilon}\,d\lambda(s)
\gtrsim\lambda(I(t,r)).
\]
The square estimate \eqref{eq:liftprofile}, comparability of squares and balls, and doubling to move the center by $|y|\le f(t)\le r$ show that
\[
\sigma(B((t,y),r))\gtrsim\lambda(I(t,r)).
\]
This also holds at $t=0$. The two hypotheses of Lemma~\ref{lem:whitneyreplacement} are satisfied. That lemma produces $\nu\in\D(\R^2)$ with $\nu\restr E_f=\mu$ and absolute continuity outside $E_f$.
\end{proof}

\begin{proof}[Proof of Corollary~\ref{thm:axialreduction}]
Lemma~\ref{lem:canonical} supplies the canonical doubling measure with projection $\lambda$. Theorem~\ref{thm:potentialcriterion} tests only this projection, giving the equivalence. If $c\lambda_1\le\lambda_2\le C\lambda_1$ on Borel sets, both sides of \eqref{eq:potentialtesting} change by fixed factors. The same parameters $\varepsilon,N$ therefore give equivalent tests.
\end{proof}
\subsection{Localization to finitely many cusps}

\begin{corollary}\label{cor:globalpotential}
Let $\Omega$ satisfy Definition~\ref{def:domain}, let $E=\overline\Omega$, and let $\mu\in\D(E)$. At each zero-angle cusp choose a straightened truncated model and the axial projection $\lambda_i$ of the transformed local measure. Then $\mu\in\mathscr R(E)$ if and only if each pair $(f_i,\lambda_i)$ satisfies \eqref{eq:potentialtesting} on the corresponding axial interval. If there is no zero-angle cusp, $\D(E)\subseteq\mathscr R(E)$.
\end{corollary}
\begin{proof}
Truncation at a fixed positive axial parameter preserves internal doubling. Near the new cross-section, the truncated set contains an interior ball of radius comparable to the testing radius; comparison with that ball bounds the restricted mass of the doubled ball. Elsewhere sufficiently small balls inherit the old estimate, and compactness handles the remaining scales. A global bi-Lipschitz shear preserves doubling and commutes with restriction. Thus each local measure is a doubling measure on its model.

If $\mu=\nu\restr E$ for some $\nu\in\D(\R^2)$, its transformed measure restricts to every local measure, so Theorem~\ref{thm:potentialcriterion} gives all the tests. Conversely, each test produces a reference measure $\sigma_i$ satisfying the width comparison and coarse domination in Lemma~\ref{lem:whitneyreplacement}.
Choose disjoint cusp neighborhoods, and use one Whitney decomposition of $\R^2\setminus E$. Near each cusp, use the original measure of a reflected interior ball for fine squares and $\sigma_i(F_i(Q))$ for coarse squares, where $F_i$ is the coordinate change. Elsewhere use ordinary reflected balls, and use constant density far away.

Bi-Lipschitz images of the squares have disjoint interiors, diameters and boundary distances comparable to $s_Q$, contain balls of radius comparable to $s_Q$, and satisfy the same packing bounds. These are the properties used in the replacement proof. Pullbacks of reflected balls retain comparable boundary depth, so their overlap remains bounded. The proof of that lemma therefore gives local finiteness, adjacent comparisons and doubling in each cusp chart.

At a fixed interface away from a tip, the width is bounded below. Every sufficiently small square there uses reflection on both sides, and internal doubling compares any two choices of reflected balls. Only finitely many positive-size interface squares remain. A finite cover of the rest of the boundary supplies the estimates in Lemma~\ref{lem:corkrealization}. The same compact-scale and far-field argument gives global doubling. All new mass is placed outside $E$, so $\nu\restr E=\mu$. If there are no zero-angle cusps, the interior ball condition noted after Definition~\ref{def:domain} and Lemma~\ref{lem:corkrealization} apply directly.
\end{proof}

\subsection{Doubling widths and power-comparable cusps}
\begin{proposition}\label{prop:widthtest}
Let $f\in C^1([0,a])$ satisfy $f(0)=f'(0)=0$ and $f(t)>0$ for $t>0$.
Suppose that, on some interval $[0,a_0]$, the function $f$ is nondecreasing and
\[
f(2t)\le Kf(t)\qquad(0<t<a_0/2)
\]
for a fixed $K\ge1$. Every doubling measure $\lambda$ on $[0,a]$ satisfies
the potential test in Theorem~\ref{thm:potentialcriterion}.
Consequently, $\D(E_f)\subseteq\mathscr R(E_f)$.
\end{proposition}
\begin{proof}
Write $D=D_\lambda$ and $d=\log_2D$. There are constants $q>0$ and $C_0$
such that
\begin{equation}\label{eq:axialreversewidth}
\lambda([0,s])\le C_0(s/t)^q\lambda([0,t]),
\qquad 0<s\le t\le a_0.
\end{equation}
Indeed, the interval $I(3R/4,R/8)$ lies in $(R/2,R)$, whereas
$[0,R]\subset I(3R/4,2R)$. Four doubling steps imply
\[
\lambda([0,R/2])\le(1-D^{-4})\lambda([0,R]).
\]
Iterating proves \eqref{eq:axialreversewidth} with
$q=-\log_2(1-D^{-4})$. In particular, $\lambda(\{0\})=0$.

Set $p=\max\{2,\log_2K\}$. Monotonicity and iteration of width doubling give
\begin{equation}\label{eq:relativewidthtest}
f(s)\ge K^{-1}(s/t)^p f(t),\qquad 0<s\le t\le a_0.
\end{equation}
Choose $0<\varepsilon<\min\{1,q/p\}$ and an integer $N>\max\{1,d\}$.
For fixed sufficiently small $t>0$, put $u=f(t)$; then $u<t/4$.
On the dyadic intervals $[2^{-j-1}t,2^{-j}t]$, equations
\eqref{eq:axialreversewidth} and \eqref{eq:relativewidthtest} yield
\[
\begin{split}
\int_{(0,t]} f(s)^{-\varepsilon}\,d\lambda(s)
&\le C u^{-\varepsilon}\lambda([0,t])
 \sum_{j=0}^{\infty}2^{-j(q-p\varepsilon)}\\
&\le C' u^{-\varepsilon}\lambda([0,t]).
\end{split}
\]
This also proves integrability of $f^{-\varepsilon}$ near zero; away from zero,
$f$ has a positive lower bound on every compact subinterval.

For completeness, the ordinary doubling-kernel estimate is
\begin{equation}\label{eq:widthkernelestimate}
\int_{[0,a]}\left(\frac{u}{u+|t-s|}\right)^N d\lambda(s)
\le C\lambda(I(t,u)).
\end{equation}
It follows by splitting into $I(t,u)$ and the annuli
$2^ju\le |t-s|<2^{j+1}u$: the resulting series is bounded by
$C\sum_{j\ge0}2^{-j(N-d)}\lambda(I(t,u))$.

On $0<s<t/2$, the kernel is at most $(2u/t)^N$. Also, doubling gives
\[
\lambda([0,t])\le C(t/u)^d\lambda(I(t,u)).
\]
The preceding integrability estimate therefore implies
\[
\begin{split}
\int_{(0,t/2)}\left(\frac{u}{u+|t-s|}\right)^N
f(s)^{-\varepsilon}\,d\lambda(s)
&\le C (u/t)^{N-d}u^{-\varepsilon}\lambda(I(t,u))\\
&\le C u^{-\varepsilon}\lambda(I(t,u)).
\end{split}
\]
For $s\ge t/2$, local monotonicity and width doubling give
$f(s)\ge f(t/2)\ge u/K$ as long as $s\le a_0$.
On $[a_0,a]$, positivity and compactness give the same estimate with a
possibly different constant, uniformly for sufficiently small $t$.
Consequently, \eqref{eq:widthkernelestimate} gives
\[
\int_{[t/2,a]}\left(\frac{u}{u+|t-s|}\right)^N
f(s)^{-\varepsilon}\,d\lambda(s)
\le C u^{-\varepsilon}\lambda(I(t,u)).
\]
Together these inequalities are exactly the required potential test.
For $\mu\in\mathcal D(E_f)$, its axial projection is doubling by Lemma~\ref{lem:projection}; applying the potential criterion completes the proof.
\end{proof}

\begin{proof}[Proof of Theorem~\ref{thm:power}]
Straighten the midpoint of the two arcs by the area-preserving shear in Lemma~\ref{lem:straight}. Apply Lemma~\ref{lem:transverse} with $g(t)=t^\alpha$. Lemma~\ref{lem:powerarea} gives local doubling of restricted area on the power cusp, and bounded area distortion transfers that property to the cusp of width $f$. Away from $0$ the boundary and its nonzero-angle junctions satisfy the corkscrew condition, so restricted area is locally doubling there as well. The localization observation at the start of Section~2 gives global doubling; \eqref{eq:openclosed} gives the open-domain assertion.

The model width $g(t)=t^\alpha$ satisfies Proposition~\ref{prop:widthtest}.
The comparable-width bi-Lipschitz map of Lemma~\ref{lem:transverse}
transfers the inclusion $\D(E_g)\subseteq\mathscr R(E_g)$ to $E_f$. Away from zero the ratio
$f/g$ has positive finite bounds by compactness, so the comparison can be
used on the whole truncated model. The localization argument of Corollary~\ref{cor:globalpotential} gives $\D(\Omega)\subseteq\mathscr R(\Omega)$ for the original domain. The strict inclusion follows from Theorem~\ref{thm:bad}.
\end{proof}

\subsection{Controlled variable powers}

\begin{theorem}\label{thm:variablepower}
Suppose $a\in C^1((0,t_0])$, where $0<t_0<1$, satisfies
\begin{equation}\label{eq:variablepowercondition}
 1<a_-\le a(t)\le a_+<\infty,\qquad |ta'(t)|\le M.
\end{equation}
Set $f(0)=0$ and $f(t)=t^{a(t)}$ for $0<t\le t_0$. Then $\D(E_f)\subseteq\mathscr R(E_f)$. For every $\mu\in\D(E_f)$, a measure $\nu\in\D(\R^2)$ with $\nu\restr E_f=\mu$ can be chosen absolutely continuous on $\R^2\setminus E_f$.
\end{theorem}

\begin{proof}
The estimate
\[
 |f'(t)|\le t^{a_--1}\big(a_++M\log(1/t)\big)\longrightarrow0
\]
shows that $f\in C^1([0,t_0])$ and $f'(0)=0$. Fix $\mu\in\D(E_f)$ and
write $\lambda=\pi_\#\mu$. By Lemma~\ref{lem:projection}, $\lambda$ is
full-support doubling on $[0,t_0]$. Such a measure has no atoms: a
separated subinterval in each sufficiently small interval has a fixed
positive fraction of the larger interval's mass, so iteration makes
singleton masses zero. In particular, there is $\delta>0$ such that
\begin{equation}\label{eq:varaxialdecay}
 \lambda([0,u])\le C(u/R)^\delta\lambda([0,R]),
 \qquad 0<u\le R\le t_0.
\end{equation}
Reflect $\lambda$ across the two endpoints and repeat with period $2t_0$,
obtaining $\widetilde\lambda\in\D(\R)$. Indeed, intervals crossing a
reflection point are compared with one-sided intervals by doubling;
intervals shorter than a period meet a bounded number of reflected
pieces, and longer intervals are controlled by the common period mass.
There is no endpoint mass to duplicate.

Choose
\[
 0<\kappa<\min\{1,a_--1,\delta\},\qquad
 \beta(t)=\frac{\kappa}{a(t)-1}.
\]
Then $0<\beta_-\le\beta(t)\le\beta_+<1$ and
$|t\beta'(t)|\le C$. Near the origin define
\begin{equation}\label{eq:varreferenceweight}
 W(t,y)=
 \begin{cases}
 t^{-1-\kappa}
 \left(\dfrac{\max\{|y|,f(t)\}}{t}\right)^{\beta(t)-1},
                  &t>0,\ |y|<t,\\[5pt]
 |y|^{-1-\kappa}, &t>0,\ |y|\ge t,\\[2pt]
 \max\{|t|,|y|\}^{-1-\kappa},&t\le0,
 \end{cases}
\end{equation}
with an arbitrary value at $(0,0)$. Continue $W$ to a positive function
$w$ equal to $1$ outside a small neighborhood of the origin. The
continuation can be chosen bounded above and below on the intervening
annulus. Put
\begin{equation}\label{eq:varreference}
 d\sigma(t,y)=w(t,y)\,d\widetilde\lambda(t)\,dy.
\end{equation}
We verify its local doubling estimates rather than assuming that a
product with an arbitrary weight preserves doubling.

Since $(f(t)/t)^{\beta(t)}=t^\kappa$, one has
\begin{equation}\label{eq:vartrace}
 W(t,y)=f(t)^{-1}\qquad (|y|\le f(t)).
\end{equation}
For $0<r\le t$, direct integration gives
\begin{equation}\label{eq:varvertical}
 \int_{-r}^{r}W(t,y)\,dy\asymp rW(t,r).
\end{equation}
For $r\le f(t)$ this is immediate. Otherwise put
$H=t^{-\kappa}(r/t)^{\beta(t)}\ge1$; the integral is
$2\{1+(H-1)/\beta(t)\}$ and $rW(t,r)=H$.

For $r\le t/4$ and $|s-t|\le r$, logarithmic differentiation gives
\begin{equation}\label{eq:varhorizontal}
 W(s,r)\asymp W(t,r).
\end{equation}
On the branch $r<f(s)$,
$r|\partial_s\log W(s,r)|=r|f'(s)|/f(s)\le |f'(s)|$.
On the branch $f(s)<r<s$,
\[
 \partial_s\log W(s,r)
   =\frac{-\kappa-\beta(s)}{s}+\beta'(s)\log(r/s).
\]
Integration over an interval of length at most $2r$ gives a bound
$C(r/t)(1+\log(t/r))\le C$ for the logarithmic oscillation.
The formulas agree at the branch interfaces, so the comparison holds
across them as well. It follows that
\begin{equation}\label{eq:varaxisprofile}
 \sigma\big((t-r,t+r)\times(-r,r)\big)
 \asymp rW(t,r)\lambda(I(t,r)),\qquad r\le t/4,
\end{equation}
whenever the square lies in the neighborhood where $w=W$.
Fixed-factor changes in $r$ change $rW(t,r)$ by bounded factors:
it is $r/f(t)$ below $f(t)$ and
$t^{-\kappa}(r/t)^{\beta(t)}$ above $f(t)$.
Together with doubling of $\lambda$, this proves doubling for
axis-centered balls whose radii are small compared with $t$.

For a ball whose radius is a sufficiently small fraction of its distance
from the horizontal axis, the same logarithmic derivative estimates
make $w$ comparable throughout the ball and its double. If $|y|\ll t$,
the relevant additional bound is
$(r/t)\log(t/|y|)\le C$ when $r\le |y|/4$;
if $|y|\asymp t$ or $t\le0$, the radial-power formula applies.
The doubling property of $d\widetilde\lambda(t)\,dy$ then gives the
required comparison. Balls whose distance from the axis is comparable
to their radius are bounded above by an axis-centered square and below
by a fixed relative subball away from the axis; the estimates just
proved and \eqref{eq:varaxisprofile} apply.

It remains to consider the origin. Column integration in
\eqref{eq:varreferenceweight}, reflection of $\lambda$, and
\eqref{eq:varaxialdecay} give
\begin{align}
 \sigma([-r,r]^2)
 &\asymp\int_{(0,r]}t^{-\kappa}\,d\lambda(t)
 \asymp r^{-\kappa}\lambda([0,r]).\label{eq:varorigin}
\end{align}
For the upper bound, split $(0,r]$ into dyadic intervals and sum
$C r^{-\kappa}\lambda([0,r])
\sum_{j\ge0}2^{-j(\delta-\kappa)}$; the lower bound is immediate.
This proves local finiteness and doubling at the origin.
If the distance of a center from the origin is at most a fixed multiple
of the radius, a containing origin-centered square gives the upper
bound. For the lower bound, use
$W(t,y)\ge\max\{|t|,|y|\}^{-1-\kappa}$ on an inscribed rectangle and
compare its axial interval with $[0,r]$ using doubling of
$\widetilde\lambda$. These arguments cover all small balls. The positive
bounded continuation and the product measure away from the origin
complete the proof that $\sigma\in\D(\R^2)$.

We now check the hypotheses of Lemma~\ref{lem:whitneyreplacement}.
Let $p=(t,y)\in E_f$, with $t+r$ sufficiently small. For every fixed
$0<c_1<c_2<\infty$,
\begin{equation}\label{eq:varinterface}
 c_1f(t)\le r\le c_2f(t)
 \quad\Longrightarrow\quad
 \sigma(B(p,r))\asymp\lambda(I(t,r))
                    \asymp\mu(B_{E_f}(p,r)).
\end{equation}
Indeed, $f(t)=o(t)$ makes \eqref{eq:varaxisprofile} applicable, and
$rW(t,r)$ is bounded above and below in this fixed-ratio range.
Doubling moves the center between $p$ and $(t,0)$; the last comparison
is the full-section strip estimate underlying Lemma~\ref{lem:projection}.
More generally, for every fixed $c>0$,
\begin{equation}\label{eq:varcoarsedominance}
 r\ge cf(t)\quad\Longrightarrow\quad
 \mu(B_{E_f}(p,r))\asymp\lambda(I(t,r))
                         \lesssim\sigma(B(p,r)).
\end{equation}
When $r\ll t$, this follows from \eqref{eq:varaxisprofile} and
$rW(t,r)\gtrsim1$; when $t\lesssim r$, it follows from
\eqref{eq:varorigin} and doubling. At $t=0$, the latter estimate applies
directly. The replacement lemma therefore yields an ambient doubling
measure equal to $\mu$ on every Borel subset of $E_f$, with area density
on its complement. Notice that \eqref{eq:vartrace} alone gives only a
reference measure and does not yet imply $\sigma\restr E_f=\mu$; the replacement
step is needed to preserve the original transverse mass distribution.
\end{proof}

\begin{example}\label{ex:oscillatorywidth}
The width
\[
 f(t)=t^{3+\sin\log(1/t)}\qquad(0<t\le t_0)
\]
satisfies Theorem~\ref{thm:variablepower}, since $2\le a(t)\le4$ and
$|ta'(t)|\le1$. Consequently $\D(E_f)\subseteq\mathscr R(E_f)$. Nevertheless,
for $t_k=e^{-2\pi k}$,
\[
 \log\frac{f(2t_k)}{f(t_k)}
 =2\pi k\sin(\log2)+(3-\sin(\log2))\log2\longrightarrow\infty.
\]
Thus a doubling condition on the width is not necessary for $\D(E_f)\subseteq\mathscr R(E_f)$.
\end{example}

\subsection{Nonempty intersection without either inclusion}
The existence of a compact set of upper box dimension less than one with positive mass for a doubling measure, and the symmetric ternary Bernoulli construction used below, are known; see Rossi--Shmerkin~\cite[Section~4]{RS}. We include the required special case to explain the endpoint gluing. The lemma itself is not claimed as a new result.

\begin{lemma}\label{lem:compacttypical}
There exist a finite doubling measure $\tau$ with full support on $[0,1]$ and a compact set
$K\subset(1/4,3/4)$ such that
\[
 \tau(K)>0,\qquad \overline{\dim}_{\rm B}K<1,
 \qquad \tau([0,r])\asymp r,\quad \tau([1-r,1])\asymp r
 \quad(0<r\le1).
\]
\end{lemma}
\begin{proof}
Let $\omega$ be the ternary Bernoulli measure with digit weights $(1/4,1/2,1/4)$.
The masses of adjacent ternary intervals of the same generation have ratio between $1/2$ and $2$: the trailing digits $2$ and $0$ involved in a carry have equal weights, so only the first changed digit contributes to the ratio.
Comparison of adjacent intervals and ternary coverings show that $\omega$ has full support and is doubling.
Its Shannon entropy is
$H=-2(1/4)\log(1/4)-(1/2)\log(1/2)<\log3$.

Choose $H<q<\log3$. By the law of large numbers for digit frequencies, for sufficiently large $n_0$ there is a compact set $K_0$ of positive $\omega$-mass such that, for every $n\ge n_0$, all length-$n$ prefixes occurring in $K_0$ have frequencies sufficiently close to $(1/4,1/2,1/4)$.
At generation $n$, the number of admissible prefixes is at most $Ce^{qn}$; any polynomial factor is absorbed by slightly increasing $q$, still keeping $q<\log3$.
Consequently, $\overline{\dim}_{\rm B}K_0\le q/\log3<1$.
By restricting to an interior closed interval and increasing $n_0$ if necessary, the compact set may be chosen inside a prescribed interior open interval while retaining positive mass.
Points with two ternary expansions form a countable null set. A compact set can be defined directly by a decreasing intersection of unions of closed basic intervals.

Set $\alpha=\log4/\log3>1$. Then
$\omega([0,r])\asymp r^\alpha$, with the analogous estimate at the right endpoint.
Affinely rescale $\omega$ to $[1/3,2/3]$, denote the resulting probability measure by $\sigma_0$, and set
$\tau=\mathcal L^1|_{[0,1]}+\sigma_0$. Its total mass is $2$.
For a ball and its double contained in $[1/3,2/3]$, doubling follows by adding the two doubling estimates.
If the doubled ball crosses an endpoint of the support of $\sigma_0$ while the smaller ball remains inside that support, the center is within a fixed multiple of the radius from the endpoint. The $\sigma_0$-mass of the doubled ball is therefore at most $Cr^\alpha\le Cr$, and is controlled by the Lebesgue part.
The other configurations crossing a support endpoint follow from the same endpoint bound. Away from the support, only Lebesgue measure remains.
Thus $\tau$ is doubling and agrees with Lebesgue measure in fixed neighborhoods of $0$ and $1$.
Apply the same affine rescaling to the compact typical set $K_0$ constructed above. The resulting set satisfies
$K\subset[1/3,2/3]\subset(1/4,3/4)$, $\tau(K)>0$, and has upper box dimension less than one.
No quasisymmetric change of coordinates is needed.
\end{proof}

\begin{proof}[Proof of Theorem~\ref{thm:partialc1}]
Fix $K,\tau$ as in the lemma, and choose $d<1$ such that
\begin{equation}\label{eq:minkowskithin}
 |\{z\in[0,1]:\operatorname{dist}(z,K)<h\}|\le C h^{1-d}
 \quad(0<h<1).
\end{equation}
Choose $a>1$. Let $d_K(z)=\operatorname{dist}(z,K)$ on the whole real line, and
let $\varphi$ be a nonnegative smooth mollifier supported on $[-1,1]$ with integral $1$. Write
$q_\eta=d_K*\varphi_\eta$.
Since the distance function is $1$-Lipschitz,
\[
 |q_\eta-d_K|\le\eta,\qquad |q'_\eta|\le1,
 \qquad (q_\eta^2+\eta^2)^{1/2}\asymp d_K+\eta.
\]
Choose a smooth cutoff $\chi$ equal to $1$ on an interior closed interval containing $K$, and equal to $0$ near the endpoints, with a fixed positive distance from $K$ wherever $\chi\ne1$.
For $0<\eta\le1$, set
\[
 G_\eta(z)=\chi(z)(q_\eta(z)^2+\eta^2)^{a/2}+1-\chi(z).
\]
Then $G_\eta\in C^\infty([0,1])$, it equals $1$ near the endpoints, and
\begin{equation}\label{eq:Gbounds}
 G_\eta(z)\asymp(\operatorname{dist}(z,K)+\eta)^a,
 \qquad \|G_\eta\|_\infty+\|G'_\eta\|_\infty\le C_a.
\end{equation}
All constants are independent of $\eta$. The derivative bound follows from
$q_\eta(q_\eta^2+\eta^2)^{a/2-1}\lesssim(q_\eta+\eta)^{a-1}$.
Moreover, $G_\eta\le C\eta^a$ on $K$.

Set $t_k=4^{-k}$, $I_k=[t_k,2t_k]$, and $\eta_k=e^{-1/t_k}$.
On $I_k$, define
\[
 f(t)=t^2G_{\eta_k}\left(\frac{t-t_k}{t_k}\right),
\]
Set $f(t)=t^2$ elsewhere, and $f(0)=0$.
Since $G_{\eta_k}=1$ near both ends of each block, $f$ is smooth at every $t>0$. Moreover,
$|f'(t)|\le Ct_k$ on $I_k$. Hence $f'(0)=0$, and $f$ is positive on $(0,1]$.

\emph{A measure in $\D(E_f)\setminus\mathscr R(E_f)$.}
On each $I_k$, use the affine copy of $\tau$ multiplied by $t_k$, and use Lebesgue measure outside $\bigcup_k I_k$. Denote the resulting measure by $\lambda_b$. Each complete block has mass $2t_k$. Equivalently, $\lambda_b=\mathcal L^1+\sum_k t_k(A_k)_\#\sigma_0$, where $A_k(z)=t_k+t_kz$.
This measure is doubling, with a constant independent of the block index.
Indeed, adjacent intervals contained in a single block are compared using doubling of $\tau$. Intervals touching a block endpoint have mass comparable to their length by
$\tau([0,r])\asymp r$ and its right-endpoint analogue. An interval crossing several blocks is decomposed into complete blocks and at most two end pieces, yielding the same comparison.
If a small interval lies in one block but a fixed dilate exits the block, first use doubling inside the block to compare with an interval at the endpoint, then apply the endpoint estimate.
These cases cover every ball and its double.

Define
\[
 d\mu_b(t,y)=\frac{1}{2f(t)}\,d\lambda_b(t)\,dy\qquad(|y|\le f(t)).
\]
The standard ball estimate for a cusp with small slope gives
\[
\mu_b(B_{E_f}((t,y),r))\asymp\lambda_b(I(t,r))\min\{1,r/f(t)\},
\]
so $\mu_b$ is doubling.
Write $K_k=t_k+t_kK$. Then $\lambda_b(K_k)=t_k\tau(K)$ and
$f(t)\le C t_k^2\eta_k^a$ on $K_k$.
For every $\varepsilon>0$,
\[
 \int f(t)^{-\varepsilon}\,d\lambda_b(t)
 \ge c\sum_k t_k(t_k^2\eta_k^a)^{-\varepsilon}
 =c\sum_k t_k^{1-2\varepsilon}e^{a\varepsilon/t_k}=\infty.
\]
Corollary~\ref{thm:internalnonrestriction} forces the axial projection of every measure in $\D(E_f)\cap\mathscr R(E_f)$ to have some finite negative width moment. Thus $\mu_b\in\D(E_f)\setminus\mathscr R(E_f)$.

\emph{A measure in $\D(E_f)\cap\mathscr R(E_f)$.}
Take $\lambda_g=dt$ and
$0<\varepsilon<\min\{1/2,(1-d)/a\}$.
The layer-cake formula applied to \eqref{eq:minkowskithin} gives
\[
 \sup_{0<\eta\le1}\int_0^1G_\eta(z)^{-\varepsilon}\,dz<\infty.
\]
Thus $\rho=f^{-\varepsilon}dt$ is finite, since
$\rho(I_k)\le C t_k^{1-2\varepsilon}$ and $f=t^2$ off the blocks.
We verify the potential test. Choose
\[
 N>\max\{2,\ a(1-\varepsilon)/(a-1)\}.
\]
Fix $\theta>0$ sufficiently small that the $2\theta$-neighborhood of $K$ is compactly contained in $(0,1)$. If $t\in I_k$ and its normalized coordinate lies within distance $\theta$ of $K$, set
$\delta=\operatorname{dist}((t-t_k)/t_k,K)+\eta_k$.
Then $r=f(t)\asymp t_k^2\delta^a$.
For $|s-t|<c t_k\delta$, we have $f(s)\asymp f(t)$, so
\[
 \int_{|s-t|<ct_k\delta}\left(\frac r{r+|s-t|}\right)^N d\rho(s)
 \le C f(t)^{-\varepsilon}r=C r^{1-\varepsilon}.
\]
On the remaining part of the same block, the bound $\rho(I_k)\le C t_k^{1-2\varepsilon}$ gives
\[
 \int_{I_k\setminus(t-ct_k\delta,t+ct_k\delta)}
     \left(\frac r{r+|s-t|}\right)^N d\rho(s)
 \le C t_k^{N+1-2\varepsilon}\delta^{(a-1)N}
 \le C t_k^{2-2\varepsilon}\delta^{a(1-\varepsilon)}
 \asymp C r^{1-\varepsilon}.
\]
Here we used $N>1$ and $(a-1)N\ge a(1-\varepsilon)$.
Outside the block, the distance from $t$ is at least $c t_k$, so the contribution is at most
\[
 C(r/t_k)^N\rho([0,1])
 \le C t_k^N\delta^{aN}\le C r^{1-\varepsilon},
\]
because $N\ge2-2\varepsilon$ and $N\ge1-\varepsilon$.

If $t$ is outside the normalized $\theta$-neighborhood described above, then $f(t)\asymp t^2$. There is a uniform $c>0$ such that $f(s)\asymp t^2$ whenever $|s-t|<ct$. Inside a block, choose $c$ smaller than a fixed multiple of $\theta$; near a block endpoint, use the fixed collar on which $G_\eta=1$. The gaps between blocks have size comparable to the corresponding $t_k$, so the same choice covers points outside the blocks.
The near contribution is controlled by the same kernel integral. For the far contribution, use
$(r/t)^N\rho([0,1])\lesssim r^{1-\varepsilon}$.
Consequently, for every sufficiently small $t>0$,
\[
 \int\left(\frac{f(t)}{f(t)+|t-s|}\right)^Nf(s)^{-\varepsilon}\,ds
 \le C f(t)^{1-\varepsilon}
 \asymp C f(t)^{-\varepsilon}|I(t,f(t))|.
\]
The axial potential criterion now shows that
$d\mu_g=(2f(t))^{-1}dt\,dy$ belongs to $\D(E_f)\cap\mathscr R(E_f)$.

Finally, the graphs $y=\pm f(t)$ are regular $C^1$ arcs. Joining them by the vertical segment at $t=1$ gives a Jordan boundary consisting of finitely many regular $C^1$ arcs.
The example therefore satisfies Definition~\ref{def:domain}. Theorem~\ref{thm:bad} also gives $\mathscr R(E_f)\setminus\D(E_f)\ne\varnothing$. The open--closed correspondence yields the same three nonempty sets on $\Omega_f$.
\end{proof}

This construction is consistent with the dichotomy for the smooth subclass. Indeed, at the gap points $t=3t_k$ we have $f(t)/t^2=1$, whereas, for a fixed $z\in K$ and $t=t_k(1+z)$, we have $f(t)/t^2\le C\eta_k^a\to0$. Thus $f(t)/t^2$ has no limit at the origin, and $f\notin C^2([0,1])$. Since the horizontal coordinate of the graph has nonzero derivative, no regular reparametrization can turn this boundary into a regular $C^\infty$ arc at the endpoint.

\subsection{Smooth contact and finitely many cusps}

\begin{proof}[Proof of Theorem~\ref{thm:smooth}]
If $m<\infty$, Taylor's theorem and positivity give $f(t)=bt^m+o(t^m)$ with $b>0$. Theorem~\ref{thm:power} applies. If $m=\infty$, Taylor's theorem at every finite order gives \eqref{eq:flat}, and Theorem~\ref{thm:flat} applies.
\end{proof}

Smoothness here includes all one-sided derivatives at the cusp endpoint. The weaker condition $C^1([0,a])\cap C^\infty((0,a])$ holds for both examples in Theorems~\ref{thm:counter} and~\ref{thm:partialc1}. Regular smooth changes of coordinates preserve finite versus infinite contact order: the graph gap changes by a nonvanishing smooth factor and an invertible tangential reparametrization.

\begin{remark}\label{rem:finite}
For a piecewise $C^1$ Jordan domain with finitely many zero-angle cusps, no such cusp gives $\mathscr R(\Omega)=\D(\Omega)$. If there is at least one and every cusp width is power-comparable, then $\D(\Omega)\subsetneq\mathscr R(\Omega)$ and restricted area is doubling. A single infinitely flat cusp, or one with the obstructing geometry of Theorem~\ref{thm:counter}, gives $\mathscr R(\Omega)\cap\D(\Omega)=\varnothing$.

For finitely many regular $C^\infty$ arcs smooth up to their endpoints, this is a complete geometric classification: no zero-angle cusp gives equality; at least one cusp, all of finite contact order, gives strict inclusion $\D(\Omega)\subsetneq\mathscr R(\Omega)$; an infinite-order cusp gives disjointness. Indeed, apply Taylor's theorem at each cusp, Corollary~\ref{cor:globalpotential} and Theorem~\ref{thm:power} for the positive inclusion, and the local constructions of Theorems~\ref{thm:bad}, \ref{thm:flat}, and~\ref{thm:counter} for the remaining assertions. A finite cover handles restricted area.
\end{remark}

\begin{corollary}\label{cor:analyticjordan}
For a Jordan domain bounded by finitely many regular $C^\infty$ arcs smooth up to their endpoints, the following are equivalent:
\[
\mathscr R(\Omega)\cap\D(\Omega)\ne\varnothing,
\qquad \D(\Omega)\subseteq\mathscr R(\Omega),
\qquad\text{all zero-angle cusps have finite contact order}.
\]
In particular, $\D(\Omega)\subseteq\mathscr R(\Omega)$ if the boundary consists of finitely many regular real-analytic arcs defined on neighborhoods of their closed parameter intervals.
\end{corollary}
\begin{proof}
The equivalence follows from Remark~\ref{rem:finite} and $\D(\Omega)\ne\varnothing$. A positive analytic gap between distinct incident graphs cannot vanish to infinite order, since otherwise it vanishes identically near the endpoint.
\end{proof}

\subsection{A discrete criterion and finite certificates}\label{sec:discrete}

No monotonicity, doubling condition, or leading Taylor term is assumed for the cusp widths in this section. Let $\Omega$ be a piecewise $C^1$ Jordan domain, and let $\mu\in\D(\Omega)$. We give a necessary and sufficient condition involving only grid masses. The condition depends on the prescribed measure $\mu$; it does not constitute a classification of domains solely by boundary geometry.

\subsubsection{The grid and finite tests}

Let $\mathscr G$ be the standard half-open dyadic squares in the plane:
\[
 Q_{j,k}=2^{-j}(k+[0,1)^2),\qquad j\in\mathbb Z,\ k\in\mathbb Z^2.
\]
Write $j(Q)=j$. Each square has four children. Two squares at the same level are called adjacent if their closures share an edge. Fix a square with $\overline{Q_0}\subset\Omega$, and write $M_0=\mu(Q_0)>0$.

For $C\ge1$, set $A_C=4C^2$. For each $Q$, let
$j_* =\min\{j(Q),j(Q_0)\}$, and let $P,P_0$ be the ancestors of $Q,Q_0$, respectively, at level $j_*$. Let $d_Q$ be the Manhattan distance between the integer grid indices of $P$ and $P_0$, and set
\begin{equation}\label{eq:dyadicbounds}
 a_Q(C)=C^{-d_Q}A_C^{-(j(Q)-j_*)},\qquad
 b_Q(C)=C^{d_Q}A_C^{j(Q_0)-j_*}.
\end{equation}
These quantities are determined explicitly by the grid and $C$, and satisfy $0<a_Q(C)\le b_Q(C)<\infty$.

Introduce real unknowns $m_Q$, $Q\in\mathscr G$, and consider the countable linear system $\mathcal L_C(\mu)$:
\begin{align}
 m_Q&=\sum_{Q'\text{ a child of }Q}m_{Q'},
                       &&Q\in\mathscr G,\label{eq:dyadicconservation}\\
 C^{-1}m_{Q'}&\le m_Q\le Cm_{Q'},
                       &&Q,Q'\text{ adjacent at the same level},\label{eq:dyadicneighbor}\\
 m_Q&=\mu(Q),           &&\overline Q\subset\Omega,\label{eq:dyadicdata}\\
 a_Q(C)M_0&\le m_Q\le b_Q(C)M_0,
                       &&Q\in\mathscr G.\label{eq:dyadicbox}
\end{align}
We say that every finite subsystem is feasible if each finite selection of the displayed equations and inequalities has a real solution. For fixed $C$, this is an ordinary finite-dimensional linear feasibility problem. Solutions to different finite subsystems need not be consistent in advance.

\begin{theorem}\label{thm:dyadiccriterion}
The following statements are equivalent for $\mu\in\D(\Omega)$:
\begin{enumerate}
\item $\mu\in\mathscr R(\Omega)$.
\item There exists a constant $C\ge1$ for which the entire system $\mathcal L_C(\mu)$ is feasible.
\item There exists a constant $C\ge1$ for which every finite subsystem of $\mathcal L_C(\mu)$ is feasible.
\end{enumerate}
For necessity, one may take $C=D_\nu^4$. For sufficiency, the doubling constant of the resulting planar measure is controlled by $C$ and absolute constants.
\end{theorem}

\subsubsection{Necessity and compactness}

\begin{proof}[Proof of Theorem~\ref{thm:dyadiccriterion}: $(1)\Rightarrow(2)$]
Choose $\nu\in\D(\R^2)$ with $\nu\restr\Omega=\mu$, and set $m_Q=\nu(Q)$. An ambient doubling measure gives zero mass to lines, so the half-open convention does not affect square masses. Finite additivity gives \eqref{eq:dyadicconservation}, and $\nu\restr\Omega=\mu$ gives \eqref{eq:dyadicdata}. Inscribed balls in two squares sharing an edge, together with a fixed dilate containing both squares, give \eqref{eq:dyadicneighbor} by doubling; for example, $C=D_\nu^4$ suffices.

We must check that \eqref{eq:dyadicbox} imposes no additional restriction. Any two children of a parent square are connected by at most two edge-adjacency steps. Hence each child has mass at least $A_C^{-1}$ times the parent mass and at most the parent mass. Also, $P_0$ and $P$ are connected by a path of $d_Q$ edge-adjacency steps, so their mass ratio lies between $C^{-d_Q}$ and $C^{d_Q}$. Combining these comparisons along the ancestor chains and using $m_{Q_0}=M_0$ gives exactly \eqref{eq:dyadicbox}.
\end{proof}

Statement~$(2)$ clearly implies $(3)$. Conversely, consider the compact product space
\[
 \mathcal K_C=\prod_{Q\in\mathscr G}[a_Q(C)M_0,b_Q(C)M_0].
\]
Each conservation equation, adjacency comparison, or interior-data constraint defines a closed subset of this space and involves only finitely many coordinates. Statement~$(3)$, together with the interval bounds for the relevant coordinates, gives the finite intersection property. The total intersection is therefore nonempty, proving $(3)\Rightarrow(2)$. Since the index set is countable, the same argument can also be given by a diagonal subsequence.

One must first fix the same $C$ and then test all finite subsystems. Allowing $C$ to grow with the size of the test does not imply $\mu\in\mathscr R(\Omega)$.

\subsubsection{Constructing a planar measure}

\begin{proof}[Proof of Theorem~\ref{thm:dyadiccriterion}: $(2)\Rightarrow(1)$]
Fix a solution $(m_Q)$. In each unit square, use conservation on its four-branch tree to construct a probability measure, and multiply by the mass of that unit square. More explicitly, the conditional probability of choosing a child is $m_{Q'}/m_Q$. An infinite nested sequence of closed squares has diameters tending to zero and intersects in a unique point. Push the tree measure forward to that point.

We must exclude duplicate mass on dyadic boundaries. Within any square $Q$, each child has mass at least $A_C^{-1}m_Q$. Retaining the two children along one edge discards at least two children, so the retained mass is at most
$(1-2/A_C)m_Q$. After $k$ repetitions, the coding mass of the strip of width $2^{-k}\ell(Q)$ along that edge is at most
\[
 (1-2/A_C)^k m_Q\longrightarrow0.
\]
Thus square boundaries and all dyadic grid lines are null sets. The construction gives a Borel measure satisfying $\nu(Q)=m_Q$. Summing over the unit squares gives a measure on the whole plane. Repeated conservation also gives the prescribed masses of squares at negative levels.

Every compact set meets only finitely many unit squares, so $\nu$ is locally finite. A locally finite Borel measure on Euclidean space is Radon. Every open ball contains a dyadic square, whose mass is positive; hence $\nu$ has full support.

To verify doubling, fix $B(x,r)$ and choose a level with $r/8\le\ell(Q)<r/4$. Let $Q_x$ be the square at that level containing $x$. Then $Q_x\subset B(x,r)$. At most $N_0$ squares of this level cover $B(x,2r)$, and each can be joined to $Q_x$ by at most $N_1$ edge-adjacency steps, where $N_0,N_1$ are absolute constants. Thus
\[
 \nu(B(x,2r))\le N_0 C^{N_1}m_{Q_x}
 \le N_0 C^{N_1}\nu(B(x,r)).
\]
It follows that $\nu\in\D(\R^2)$.

Finally, we verify $\nu\restr\Omega=\mu$. The measure $\mu$ also gives zero mass to lines inside $\Omega$: in small balls compactly contained in $\Omega$, relative porosity of a line and doubling of $\mu$ give zero mass by the differentiation theorem. Interior dyadic squares satisfy \eqref{eq:dyadicdata}; their subdivisions generate the Borel sets in $\Omega$, and the omitted grid lines are null for both measures. Therefore $\nu\restr\Omega=\mu$, as required. The boundary argument following Definition~\ref{def:domain} gives the same conclusion with $E=\overline\Omega$ in place of $\Omega$.
\end{proof}

\begin{corollary}\label{cor:finitecertificate}
For $\mu\in\D(\Omega)$, one has $\mu\notin\mathscr R(\Omega)$ if and only if, for every integer $C\ge1$, the system $\mathcal L_C(\mu)$ has a finite infeasible subsystem.
\end{corollary}
\begin{proof}
Theorem~\ref{thm:dyadiccriterion} gives the assertion with all real $C\ge1$. If $\mu\in\mathscr R(\Omega)$, choose a planar doubling measure whose restriction is $\mu$. Any sufficiently large integer is an admissible adjacency-comparison constant, and the necessity argument gives a feasible system for that integer. This proves the stated version.
\end{proof}

Failure of a finite test for one fixed $C$ excludes only that adjacency-comparison constant; it does not establish $\mu\notin\mathscr R(\Omega)$. Likewise, passing all tests up to a finite level does not replace passing every finite subsystem. Deciding whether $\D(\Omega)\subseteq\mathscr R(\Omega)$ still requires applying the criterion to every $\mu\in\D(\Omega)$. This criterion does not by itself give a classification based solely on the geometry of the boundary.

For an arbitrary positive $C^1$ width, the potential and finite-system criteria determine membership of each prescribed internal measure in $\mathscr R(\Omega)$. They do not eliminate the measure quantifier to give a purely geometric characterization of all domains with $\mathscr R(\Omega)\cap\D(\Omega)=\varnothing$, with $\D(\Omega)\subseteq\mathscr R(\Omega)$, or with nonempty intersection and neither inclusion. In particular, the local condition~\eqref{eq:localnecessary} and the negative moment condition~\eqref{eq:negativewidthmoment} are necessary conditions only.

\begingroup\small\linespread{1}\selectfont

\endgroup
\end{document}